\documentclass[11pt]{amsart}
\usepackage{bbm}
\usepackage{mathbbold}
\usepackage{pifont}
\usepackage{mathrsfs}
\usepackage{amsfonts}
\usepackage{latexsym, amssymb, amsmath, amscd, amsfonts ,texdraw, mathrsfs, amsthm}
\usepackage{colordvi, ifpdf}
\renewcommand{\proof}{\noindent{\it Proof.\ \ }}
\def\bfl{{\bf l}}\def\bfv{{\bf v}}
\def\div{\,\,\big|\,\,}
\def\mfs{\mathfrak{s}}\def\mfl{\mathfrak{l}}\def\mfr{\mathfrak{r}}\def\mfw{\mathfrak{w}}
\def\mfD{\mathfrak{D}}
\def\mfC{\mathfrak{C}}
\def\mcS{\mathcal{S}}\def\mcB{\mathcal{B}}\def\mcP{\mathcal{P}}\def\mcE{\mathcal{E}}
\def\mfu{\mathfrak{u}}\def\mfv{\mathfrak{v}}

\def\mbk{\mathbbm{k}}\def\mbv{\mathbbm{v}}
\def\mbD{\mathbb{D}}\def\mcD{\mathcal{D}}

\newtheorem{thm}{Theorem}[section]
\newtheorem{lem}[thm]{Lemma}

\newtheorem{cor}[thm]{Corollary}
\newtheorem{prop}[thm]{Proposition}

\newtheorem{con}[thm]{Construction}
\newtheorem{ex}[thm]{Example}

\theoremstyle{definition}

\def\qed{\hfill \rule{4pt}{7pt}}

\begin{document}

\title[symmetric graphs]{A class of symmetric graphs with $2$-arc-transitive
quotients}
\thanks{This work is supported by  the 973 Program and the NSF of China}
\thanks{{\em Corresponding email}: lu@nankai.edu.cn}
\author{Bin Jia}
\address{Center for Combinatorics\\
 LPMC, Nankai University\\
 Tianjin 300071\\
P. R. China} \email{jiabinqq@gmail.com}

\author{Zai Ping Lu}
\address{Center for Combinatorics\\
 LPMC, Nankai University\\
 Tianjin 300071\\
P. R. China} \email{lu@nankai.edu.cn}

\author{Gai Xia Wang}
\address{Center for Combinatorics\\
LPMC, Nankai University\\
Tianjin 300071\\
P. R. China} \email{wgx075@163.com}

\begin{abstract}
Let $\Gamma$ be a finite $X$-symmetric graph with a nontrivial $X$-invariant partition $\mathcal {B}$ on
$V(\Gamma)$ such that $\Gamma_{\mathcal {B}}$ is a connected $(X,2)$-arc-transitive graph and $\Gamma$ is not a
multicover of $\Gamma_{\mathcal {B}}$. A characterization of $(\Gamma, X, \mathcal {B})$ was given in
\cite{OACOFSG} for the  case where  $|\Gamma(C) \cap B| = 2$ for $B\in \mathcal {B}$ and $C \in \Gamma_{\mathcal
{B}}(B)$. This motivates us   to investigate the case where $|\Gamma(C) \cap B| = 3$, that is, $\Gamma[B, C]$ is
isomorphic to one of $3K_2$, $K_{3,3} -3K_2$ and $K_{3,3}$. This investigation requires a study on
$(X,2)$-arc-transitive graphs of valency $4$ or  $7$. Based on the results in \cite{FSGTTQ(II)}, we give a
characterization of tetravalent $(X,2)$-arc-transitive graphs; and as a byproduct, we prove that every
tetravalent $(X,2)$-transitive graph is either the complete graph on $5$ vertices or a near $n$-gonal graph for
some $n\ge 4$. We show that a heptavalent $(X,2)$-arc-transitive graph $\Sigma$ can occur as $\Gamma_{\mathcal
{B}}$  if and only if $X_\tau^{\Sigma(\tau)}\cong PSL(3,2)$ for $\tau\in V(\Sigma)$.

\vskip 5pt

\end{abstract}

\maketitle

\setlength{\parskip}{5pt}
\section{Introduction}
In this paper, all graphs  are assumed to be finite, nonempty,
simple and undirected. This paper involves graphs, permutation
groups and designs, the reader is referred to \cite{GTWA}, \cite{PG}
and \cite{DT} respectively for the notation and terminology not
mentioned here.

Let $\Sigma$ be a regular graph with vertex set $V(\Sigma)$ and edge
set $E(\Sigma)$. By $val(\Sigma)$ we denote the valency of $\Sigma$.
For an integer $s \geq 1$ and an $(s+1)$-sequence $\bbalpha =
(\alpha_0, \alpha_1, \ldots, \alpha_s)$ of $V(\Sigma)$, set
$\bbalpha^{-1} := (\alpha_s, \alpha_{s-1}, \ldots, \alpha_0)$,
$\bbalpha$ is called an {\em $s$-arc} of $\Sigma$ if $\{\alpha_i,
\alpha_{i+1}\} \in E(\Sigma)$ for $i = 0, 1, \ldots, s-1$, and
$\alpha_{i-1} \neq \alpha_{i+1}$ for $i = 1, 2, \ldots, s-1$. An
$s$-arc $\bbalpha = (\alpha_0, \alpha_1, \ldots, \alpha_s)$ is
called an {\em $s$-dipath} if $\alpha_i \neq \alpha_j$ for $i,j \in
\{0, 1, \ldots, s\}$ with $i \neq j$. Evidently, $\bbalpha$ is an
$s$-arc ($s$-dipath, respectively) of $\Sigma$ if and only if
$\bbalpha^{-1}$ is an $s$-arc ($s$-dipath, respectively) of
$\Sigma$. For any $s$-dipath $\bbalpha = (\alpha_0, \alpha_1,
\ldots, \alpha_s)$ of $\Sigma$, identifying $\bbalpha$ and
$\bbalpha^{-1}$ gives rise to an {\em $s$-path} $[\alpha_0,
\alpha_1, \ldots, \alpha_s]$ of $\Sigma$. Denote by $Arc_s(\Sigma)$
($Path_s(\Sigma)$, respectively) the set of $s$-arcs ($s$-paths) of
$\Sigma$. In the case where $s = 1$, we use {\em arc} and
$Arc(\Sigma)$ in place of $1$-arc and $Arc_1(\Sigma)$, respectively.

Let $X$ be a group acting on $V(\Sigma)$. The induced action of $X$
on $V(\Sigma) \times V(\Sigma)$ is defined by $(\tau, \sigma)^x =
(\tau^x, \sigma^x)$ for $(\tau, \sigma) \in V(\Sigma) \times
V(\Sigma)$ and $x \in X$. We say that $X$ preserves the adjacency of
$\Sigma$ if $Arc(\Sigma)^x = Arc(\Sigma)$, for all $x \in X$. The
graph $\Sigma$ is said to be {\em $X$-vertex-transitive} if $X$
preserves the adjacency of $\Sigma$ and acts transitively on
$V(\Sigma)$; and  $\Sigma$ is said to be {\em
$(X,s)$-arc-transitive} ({\em $(X,s)$-arc-regular}, respectively) if
in addition the induced action of $X$ on $Arc_s(\Sigma)$ is
transitive (regular, respectively). Further, $\Sigma$ is said to be
{\em $(X,s)$-transitive} if $\Sigma$ is $(X,s)$-arc-transitive but
is not $(X,s+1)$-arc-transitive. An $(X,1)$-arc-transitive graph is
usually called an {\em $X$-symmetric graph}. For $\tau\in
V(\Sigma)$, we denote by $X\tau$ the point-stabilizer of $\tau$ in
$X$. It is well-known that, for $s \in \{1, 2\}$, an
$X$-vertex-transitive graph $\Sigma$ is $(X,s)$-arc-transitive if
and only if $X_\tau$ is $s$-transitive on the neighborhood
$\Sigma(\tau) := \{\sigma \in V(\Sigma)\div (\tau, \sigma)\in
Arc(\Sigma)\}$ of $\tau$ in $\Sigma$. The reader is referred to
\cite{AGT} for basic results about symmetric graphs.

Let $\Gamma$ be a finite $X$-symmetric graph admits a nontrivial
$X$-invariant partition $\mathcal {B}$ on $V(\Gamma)$, that is, $1 <
|B| < V(\Gamma)$ and $B^x := \{\mfv^x \div \mfv\in B\} \in \mathcal
{B}$ for  $B \in \mathcal {B}$ and $x \in X$. (Such a graph is said
to be an {\em imprimitive} $X$-symmetric graph.) The {\em quotient
graph} $\Gamma_{\mathcal {B}}$ of $\Gamma$ with respect to $\mathcal
{B}$ is defined to be the graph with vertex set $\mathcal {B}$  such
that, for $B, C \in \mathcal {B}$, $B$ is adjacent to $C$ in
$\Gamma_{\mathcal {B}}$ if and only if there exists some $\mfv \in
B$ adjacent to some $\mfu \in C$ in $\Gamma$. It is easy to see that
$X$ acts transitively on the vertex set and on the arc set of
$\Gamma_{\mathcal{B}}$, that is, $\Gamma_{\mathcal {B}}$ is
$X$-symmetric. We always assume that $\Gamma_{\mathcal {B}}$ has at
least one edge, which implies that each $B \in \mathcal {B}$ is an
independent set of $V(\Gamma)$.

It has been observed in the literature that the quotient graphs of
an $(X,2)$-arc-transitive graph are usually not
$(X,2)$-arc-transitive, and that an $X$-symmetric graph with an
$(X,2)$-arc-transitive quotient itself is not necessarily
$(X,2)$-arc-transitive. (For example, several examples are given in
\cite{FP1,FP2} for the first situation; and for the second
situation, it is shown in \cite{FSGTTQ(II)} that every connected
$(X,3)$-arc-transitive graph is a quotient graph of at least one
$X$-symmetric graph which is not $(X,2)$-arc-transitive.) This
observation gave rise to a series of intensively studies of the
following two questions (Q1) and (Q2) \cite{{OACOFSG},{FSGWTTQ}} by
investigating `local' structures of imprimitive symmetric graphs and
their quotient graphs.
\begin{itemize}
\item[(Q1)] {\rm When can $\Gamma_{\mathcal {B}}$ be
$(X,2)$-arc-transitive?}
\item[(Q2)] {\rm What information of the structure of  $\Gamma$ can we
obtain from an $(X,2)$-arc-transitive quotient $\Gamma_{\mathcal
{B}}$ of $\Gamma$?}
\end{itemize}

For $B \in \mathcal {B}$ and $\mfv \in V(\Gamma)$, we set $\Gamma(B)
:= \bigcup_{\mfu \in B}\Gamma(\mfu)$,  $\Gamma_{\mathcal {B}}(B) :=
\{C \in \mathcal {B}\div (B, C) \in Arc(\Gamma_{\mathcal {B}})\}$
and $\Gamma_{\mathcal {B}}(\mfv) := \{C\in \mathcal {B}\div \mfv \in
\Gamma(C)\}$. Let $\mathcal {D}(B) := (B, \Gamma_{\mathcal {B}}(B),
|)$ denote the incidence structure such that $\mfv | C$ for $\mfv
\in B$, $C \in \Gamma_{\mathcal {B}}(B)$ if and only if $C \in
\Gamma_{\mathcal {B}}(\mfv)$. For any $B \in \mathcal {B}$, $C\in
\Gamma_{\mathcal {B}}(B)$ and $\mfv \in B$, as $\Gamma$ is
$X$-symmetric,  $v := |B|$, $k := |\Gamma(C) \cap B|$, $r :=
|\Gamma_{\mathcal {B}}(\mfv)|$ and $b := val(\Gamma_{\mathcal {B}})$
are independent of the choice of $B$ and $\mfv$, and $\mathcal
{D}(B)$ is an $X_B$-flag-transitive $1$-$(v, k, r)$ design with $b$
blocks \cite[Lemma 2.1]{FSGTTQ(II)}. $\Gamma$ is said to be a
multicover of $\Gamma_{\mathcal {B}}$ if $k = v$. For $(B, C) \in
Arc(\Gamma_{\mathcal {B}})$, denote by $\Gamma[B,C]$ the bipartite
subgraph of $\Gamma$ induced by
 $(\Gamma(C)\cap B) \cup (\Gamma(B)\cap C)$. Then $\Gamma[B,C]$ is
 independent of the choice of $(B,C)\in Arc(\Gamma_{\mathcal {B}})$ up to isomorphism, and $X_B\cap
 X_C$ acts transitively on the edges of  $\Gamma[B,C]$.

Without doubt, the triple $(\Gamma_{\mathcal {B}}, \Gamma[B,C],
\mathcal {D}(B))$ mirrors `global' and `local' information of the
structure of $\Gamma$, which allows us to reconstruct $\Gamma$ in
some sense.
 This approach to imprimitive symmetric graphs have
received attention in the literature. Gardiner and
Praeger~\cite{Gardiner-Praeger95} suggested to analyse these three
configurations, $(\Gamma, \Gamma[B,C],\mcD(B))$, and they discussed
the case when $\Gamma$ is {\em $X$-locally primitive}, that is, the
stabilizer of a vertex $\mfv\in V(\Gamma)$ in $X$  acts primitively
on the neighbourhood $\Gamma(\mfv)$ of the vertex in $\Gamma$.
In~\cite{Gardiner-Praeger98,SGWCQ} they considered the case when the
quotient $\Gamma_\mcB$ is a complete graph and the setwise
stabiliser $X_B$ (the subgroup of $X$ fixing $B$ setwise) is
$2$-transitive on $B$.  Li, Praeger and Zhou  \cite{ACOFSGWTTQ}
proved that, if $k=v-1\ge 2$, then $\mcD(B)$ contains no repeated
blocks (that is, the subsets of $B$ incident with distinct blocks of
$\mcD(B)$ are distinct) if and only if $\Gamma_\mcB$ is $(X,
2)$-arc-transitive, and further they found an elegant construction
(called the {\em $3$-arc graph} construction) for constructing all
such graphs  from $\Gamma_\mcB$.  Iranmanesh, Praeger and Zhou
\cite{FSGWTTQ}, Lu and Zhou~\cite{FSGTTQ(II)} studied the case where
the quotient $\Gamma_\mcB$ is $(X, 2)$-arc-transitive and obtained a
series of interesting results. In particular, Lu and Zhou
~\cite{FSGTTQ(II)} found the ¡®second type¡¯ 3-arc graph
construction, which led to a classification \cite{Zhoudm} of a
family of symmetric graphs. The  reader is referred to
\cite{Zhou2002a,Zhou2002b,Zhou2004,Zhou2005,OACOFSG} for further
more developments in this topic.


In answering the above two questions, a relatively explicit
classification of $(\Gamma, X, \mathcal {B})$ has been given in
\cite{OACOFSG}, when $\Gamma_{\mathcal {B}}$ is a connected
$(X,2)$-arc-transitive graph such that $2 = k \leq v-1$. This
motivated us in this paper to investigate the case where $k = 3$,
and then we give a characterization for this case.

For a group $X$ acting on a set $V$ and a subset $B$  of $V$, denote
by $X_{(B)}$ ($X_B$, respectively) the point-wise (set-wise,
respectively) stabilizer of $B$ in $X$, and by $X_B^B$ the
permutation group induced by $X_B$ on $B$. Then $X_B^B \cong X_B /
X_{(B)}$.

The following is a summary of the main results of this paper, which
 is  a sketch of our answer to (Q1) and (Q2) in the
case where $k = 3$.   More details will be given in Theorem
\ref{thm_main}.
\begin{thm}\label{thm_summary}
Let $\Gamma$ be an $X$-symmetric graph with an $X$-invariant
partition $\mcB$  on $V(\Gamma)$ such that $k=3$ and
$val(\Gamma_\mcB)\ge 2$. Let $B \in \mathcal {B}$. If $X$ is
faithful on $V(\Gamma)$ and $\Gamma$ is not a multicover of
$\Gamma_\mcB$, then $\Gamma_{\mathcal {B}}$ is
$(X,2)$-arc-transitive if and only if one of the following cases
occurs.
\begin{enumerate}
\item[{(a)}] $(v, b, r) = (4, 4, 3)$ and $X_B^B \cong A_4$ or
$S_4$;
\item[{(b)}] $(v, b, r) = (6, 4, 2)$ and
$X_B^B \cong A_4$ or $S_4$;
\item[{(c)}] $(v, b, r) = (7, 7, 3)$ and
 $X_B^B \cong PSL(3,2)$;
\item[{(d)}] $v = 3b\geq 6$, $r = 1$ and
$X_B$ acts $2$-transitively on the blocks of $\mathcal {D}(B)$.
\end{enumerate}
\end{thm}

\section{Graphs constructed from given graphs}
In this section, we aim to restate several graphs constructed from
given graphs, as well as some of their properties, which turn out to
be useful in a further characterization of $(\Gamma, X, \mathcal
{B})$ stated in Theorem \ref{thm_summary}. Hereafter, we denote by
$\mathcal {G}$ the set of triples $(\Gamma, X, \mathcal {B})$ such
that $\Gamma$ is a finite $X$-symmetric graph with a nontrivial
$X$-invariant partition $\mathcal {B}$ on $V(\Gamma)$,
$val(\Gamma_{\mathcal {B}}) \geq 2$ and $\Gamma$ is not a multicover
of $\Gamma_{\mathcal {B}}$, and by $\hat{\mathcal {G}}$ the subset
of $\mathcal {G}$ such that $\Gamma_{\mathcal {B}}$ is connected and
$X$ acts faithfully on $V(\Gamma)$, that is, $\cap_{\mfv \in
V(\Gamma)}X_{\mfv} = 1$.

The following two propositions are quoted from \cite{FSGTTQ(II)}.
\begin{prop}\label{prop_U-graph}
Let $\Sigma$ be a finite $(X,2)$-arc-transitive graph with
$val(\Sigma)\geq 2$. Let $\Delta$ be a self-paired subset of
$Arc_3(\Sigma)$, that is, $\bbalpha^{-1} \in \Delta$ whenever
$\bbalpha \in \Delta$. Define $\gimel := \gimel(\Sigma, \Delta)$ to
be the graph with vertex set $Path_2(\Sigma)$ and edge set $
\{\{[\alpha_0, \alpha_1, \alpha_2], [\alpha_1, \alpha_2,
\alpha_3]\}\div (\alpha_0, \alpha_1, \alpha_2, \alpha_3) \in
\Delta\}$. Set $P_{\tau} := \{[\tau_1, \tau, \tau_2] \in
Path_2(\Sigma)\div \tau_1, \tau_2 \in \Sigma(\tau)\}$ for $\tau \in
V(\Sigma)$,  and $\mathcal {P} := \{P_{\sigma} \div \sigma \in
V(\Sigma)\}$. If $\Delta$ is a self-paired $X$-orbit on
$Arc_3(\Sigma)$, then $(\gimel, X, \mathcal {P}) \in \mathcal {G}$
and $\Sigma \cong \gimel_{\mathcal {P}}$.
\end{prop}

The following lemma improves \cite[Theorem 4.10]{FSGTTQ(II)}.

\begin{lem}\label{lem_U-graph}
Let $(\Gamma, X, \mathcal {B}) \in \mathcal {G}$ with $b \geq 3$ and
$r=2$. Set
$$\Delta := \left\{(C, B(\mfv),
B(\mfu), D)\left|\begin{array}{llll}(\mfv, \mfu) \in
Arc(\Gamma)\\
\mfv\in B(\mfv)\in \mathcal {B}, \mfu\in B(\mfu)\in\mathcal {B}\\
C \in \Gamma_{\mathcal {B}}(\mfv) , D \in \Gamma_{\mathcal
{B}}(\mfu), C\ne B(\mfu), D\ne B(\mfv)\end{array}\right.\right\}.$$
 Suppose that  $|\Gamma(D)\cap B_0 \cap \Gamma(C)|\ne 0$ for  any $2$-path $[D,B_0,C]$ of
$\Gamma_{\mathcal {B}}$ with a given middle vertex $B_0\in \mathcal
{B}$. Then $\Gamma_{\mathcal {B}}$ is $(X,2)$-arc-transitive and
$\lambda :=|\Gamma(D) \cap B_0 \cap \Gamma(C)|$ is independent of
the choices of $[D,B_0,C]$ and $B_0$; further, $\Delta$ is a
self-paired $X$-orbit on $Arc_3(\Gamma_{\mathcal {B}})$, and either
\begin{enumerate}
\item[{(a)}] $\lambda = 1$ and $\Gamma \cong \gimel(\Gamma_{\mathcal {B}},
\Delta)$; or
\item[{(b)}] $\lambda \geq 2$ and $\Gamma$ admits a second
nontrivial $X$-invariant partition $$\mathcal{Q}:= \{\Gamma(D) \cap
B \cap \Gamma(C)\div [D,B,C] \in Path_2(\Gamma_{\mathcal {B}})\}$$
on $V(\Gamma)$, which is a proper refinement of $\mathcal {B}$ such
that $\Gamma_\mathcal{Q}\cong \gimel(\Gamma_{\mathcal {B}},
\Delta)$.
\end{enumerate}
\end{lem}
\proof  Note that $b\ge 3$. Take three distinct blocks $C, D, E\in
\Gamma_\mathcal{B}(B_0)$. Since $|\Gamma(D) \cap B_0 \cap
\Gamma(C)|\neq 0$ and $|\Gamma(E) \cap B_0 \cap \Gamma(C)|\neq 0$,
there exist $\mfv, \mfu \in \Gamma(C) \cap B_0$ with $\mfv \in
\Gamma(D)$ and $\mfu \in \Gamma(E)$. Let $\mfv', \mfu'\in C$ be such
that $(\mfv,\mfv'), (\mfu,\mfu')\in Arc(\Gamma)$. Then
$(\mfv,\mfv')^x=(\mfu,\mfu')$ for some $x\in X$ as $\Gamma$ is
$X$-symmetric. So $\mfv^x=\mfu$ and $\mfv'^x=\mfu'$, it implies
$B_0^x=B_0$ and $C^x=C$, hence $x\in X_{B_0}\cap X_{C}$. Further $C,
D^x, E\in \Gamma_\mathcal{B}(\mfu)$, it follows that $D^x=E$ as
$r:=|\Gamma_\mathcal{B}(\mfu)|=2$. Thus $X_{B_0}\cap X_{C}$ is
transitive on $\Gamma_\mathcal{B}(B_0)\setminus\{C\}$, it follows
that $X_{B_0}$ is $2$-transitive on $\Gamma_\mathcal{B}(B_0)$.
Therefore, $\Gamma_{\mathcal {B}}$ is $(X,2)$-arc-transitive. Then,
by \cite{FSGTTQ(II)}, $\lambda \geq 1$ is a constant number; and if
$\lambda = 1$, $\Delta$ is a self-paired $X$-orbit on
$Arc_3(\Gamma_{\mathcal {B}})$ and $\Gamma \cong
\gimel(\Gamma_{\mathcal {B}}, \Delta)$. In the following we  assume
$\lambda \geq 2$.

We first show $\mathcal{Q}$ is an $X$-invariant partition of
$V(\Gamma)$. Take two arbitrary $2$-paths $[D_1, B_1, C_1]$ and
$[D_2, B_2, C_2]$ of $\Gamma_{\mathcal {B}}$. Suppose that there
exists some $\mfv \in V(\Gamma)$ such that $\mfv\in (\Gamma(D_1)\cap
B_1\cap \Gamma(C_1))\cap (\Gamma(D_2)\cap B_2\cap \Gamma(C_2))$.
Then $B_1=B_2$ and $C_i,D_i \in \Gamma_{\mathcal {B}}(\mfv)$ for $i
= 1, 2$. Since $r = 2$, we have that $\{C_1, D_1\} = \{C_2, D_2\}$,
thus $[D_1, B_1, C_1]=[D_2, B_2, C_2]$. It follows that
$\mathcal{Q}$ is a partition of $V(\Gamma)$. For any $ [D, B, C] \in
Path_2(\Gamma_{\mathcal {B}})$ and $x \in X$, we have $[D,B,C]^x =
[D^x, B^x, C^x] \in Path_2(\Gamma_{\mathcal {B}})$ and so
$(\Gamma(D) \cap B \cap \Gamma(C))^x = \Gamma(D^x) \cap B^x \cap
\Gamma(C^x) \in \mathcal{Q}$. Thus $\mathcal{Q}$ is $X$-invariant.
Noting that $\Gamma$ is not a multicover of $\Gamma_\mathcal{B}$, we
know $|B|>|\Gamma(D) \cap B \cap \Gamma(C)|:=\lambda \geq 2$, so
$\mathcal{Q}$ is a proper refinement of $\mathcal {B}$. In
particular, the pair $(\mathcal{B},\mathcal{Q})$ gives an
$X$-invariant partition $\bar{\mathcal{B}}$ of
$V(\Gamma_\mathcal{Q})$.

Consider the quotient graph
$(\Gamma_\mathcal{Q})_{\bar{\mathcal{B}}}$ of $\Gamma_\mathcal{Q}$
with respect to $\bar{\mathcal{B}}$. For any $2$-path
$[\bar{D},\bar{B},\bar{C}]$ of
$(\Gamma_\mathcal{Q})_{\bar{\mathcal{B}}}$ and any $\bar{\mfv}\in
V(\Gamma_\mathcal{Q})$, we have
$|(\Gamma_\mathcal{Q})_{\bar{\mathcal{B}}}(\bar{\mfv})|=2$ and
$|\Gamma_\mathcal{Q}(\bar{D})\cap
\bar{B}\cap\Gamma_\mathcal{Q}(\bar{C})|=1$. It follows from (a) that
$\Gamma_\mathcal{Q}\cong
\gimel((\Gamma_\mathcal{Q})_{\bar{\mathcal{B}}}, \bar{\Delta})$,
where $\bar{\Delta}=\{(\bar{C}, \bar{B}(\bar{\mfv}),
\bar{B}(\bar{\mfu}),\bar{D})\div (C,B(\mfv),B(\mfu),D)\in \Delta\}$.
Moreover, it is easily shown that $\bar{\mathcal{B}}\rightarrow
\mathcal{B},\,\bar{B}\mapsto B$ is an isomorphism from
$(\Gamma_\mathcal{Q})_{\bar{\mathcal{B}}}$ to $\Gamma_\mathcal{B}$.
Therefore,  $\Gamma_\mathcal{Q}\cong
\gimel((\Gamma_\mathcal{Q})_{\bar{\mathcal{B}}}, \bar{\Delta})\cong
 \gimel(\Gamma_{\mathcal {B}}, \Delta)$.
\qed

For a finite $X$-symmetric graph $\Sigma$ with valency no less than
three, let $J(\Sigma)$ be the set of  pairs $([\tau_1, \tau,
\tau_2], [\sigma_1, \sigma, \sigma_2])$ of $2$-paths of $\Sigma$
such that $\sigma \in \Sigma(\tau) \setminus \{\tau_1, \tau_2\}$,
$\tau \in \Sigma(\sigma) \setminus \{\sigma_1, \sigma_2\}$. A subset
$\Lambda$ of $J(\Sigma)$ is said to be self-paired if $([\tau_1,
\tau, \tau_2], [\sigma_1, \sigma, \sigma_2]) \in \Lambda$ always
implies that $([\sigma_1, \sigma, \sigma_2], [\tau_1, \tau, \tau_2])
\in \Lambda$.
\begin{prop}\label{prop_H_graph}
Let $\Sigma$ be a finite $(X,2)$-arc-transitive graph with
$val(\Sigma)\geq 3$ and let $\Lambda$ be a self-paired $X$-orbit on
$J(\Sigma)$. Define a graph $\Psi := \Psi(\Sigma, \Lambda)$ with
vertex set $Path_2(\Sigma)$ such that two $2$-paths $\bbtau,
\bbsigma$ are adjacent if and only if $(\bbtau, \bbsigma) \in
\Lambda$. Then $\Psi$ is $X$-symmetric and $\mathcal {P}$ is a
nontrivial $X$-invariant partition of $V(\Psi)$ with $\Sigma \cong
\Psi_{\mathcal {P}}$, where $\mathcal {P}$ is defined as in
Proposition \ref{prop_U-graph}.
\end{prop}

We now quote a result about $3$-arc graphs  \cite{ACOFSGWTTQ}.
\begin{prop}\label{prop_3-arc_graph}
Let $\Sigma$ be a finite $(X,2)$-arc-transitive graph with
$val(\Sigma)\geq 3$ and let $\Delta$ be a self-paired $X$-orbit on
$Arc_3(\Sigma)$. The $3$-arc graph $\Xi:=\Xi(\Sigma, \Delta)$ with
respect to $\Delta$ is defined to be the graph with vertex set
$Arc(\Sigma)$ such that two arcs $(\tau, \tau_1)$ and $(\sigma,
\sigma_1)$ of $\Sigma$ are adjacent in $\Xi$ if and only if
$(\tau_1, \tau, \sigma, \sigma_1) \in \Delta$. Then $(\Xi, X,
\mathcal {A}) \in \mathcal {G}$ and $\Sigma \cong \Xi_{\mathcal
{A}}$, where $\mathcal {A} := \{A_{\tau}\div \tau \in V(\Sigma)\}$
and $A_{\tau} := \{(\tau, \sigma)|\sigma \in \Sigma(\tau)\}$ for
$\tau \in V(\Sigma)$.
\end{prop}

\begin{lem}\label{3-arc-graph}
Let $\Sigma$, $X$, $\Delta$ and $\Xi$ be as in
Proposition~\ref{prop_3-arc_graph}. Then
$r_\mathcal{A}:=|\Xi_\mathcal{A}((\tau,\tau_1))|=val(\Sigma)-1$ and
$val(\Xi)=r_\mathcal{A}\ell$, where
$(\tau,\tau_1),\,(\tau,\sigma)\in V(\Xi)=Arc(\Sigma)$ and $\ell$ is
the valency of $\Xi[A_\tau,A_\sigma]$.
\end{lem}
\proof For any arc $(\tau,\sigma)$ of $\Sigma$, there is a $3$-arc
$(\tau_1,\tau,\sigma,\sigma_1)\in \Delta$ as $X$ acts transitively
on arcs of $\Sigma$. Then $A_\tau$ and $A_\sigma$ are adjacent in
$\Xi_\mathcal{A}$. It implies that $val(\Xi)=r_\mathcal{A}\ell$. So
it suffices to show $r_\mathcal{A}=val(\Sigma)-1$. Let
$(\sigma',\sigma'_1)\in Arc(\Sigma)$. Note that $\Delta$ is
self-paired. Then $\{(\tau,\tau_1), (\sigma',\sigma'_1)\}\in E(\Xi)$
if and only if $(\tau_1,\tau,\sigma',\sigma'_1)\in \Delta$. In
particular, if $A_{\sigma'}\in \Xi_\mathcal{A}((\tau,\tau_1))$ then
$\tau_1\ne \sigma'$ and $(\tau,\sigma')\in Arc(\Sigma)$. Then
$\sigma'$, and hence $A_{\sigma'}$, has at most $val(\Sigma)-1$
choices. So $r_\mathcal{A}\le val(\Sigma)-1$. On the other hand,
since $\Sigma$ is $(X,2)$-arc-transitive, for any $\sigma'\in
\Sigma(\tau)$ with $\sigma'\ne \tau_1$, there is some $x\in X$ such
that
$(\tau_1,\tau,\sigma,\sigma_1)^x=(\tau_1,\tau,\sigma',\sigma_1^x)\in
\Delta$. It follows that $\{(\tau,\tau_1),(\sigma',\sigma_1^x)\}\in
E(\Xi)$, and so  $A_{\sigma'}\in \Xi_\mathcal{A}((\tau,\tau_1))$.
Then $r_\mathcal{A}\ge val(\Sigma)-1$. Thus $r_\mathcal{A}=
val(\Sigma)-1$. \qed


\section{Double star graphs}
If $(\Gamma, X, \mathcal {B}) \in \mathcal {G}$ such that
$\Gamma_{\mathcal {B}}$ is $(X,2)$-arc-transitive then, by
\cite{FSGTTQ(II)}, $\Gamma$ or  a quotient of $\Gamma$ is isomorphic
to one of $|E(\Gamma_\mcB)|K_2$, $\gimel(\Gamma_{\mathcal {B}},
\Delta)$, $\Psi(\Gamma_{\mathcal {B}}, \Lambda)$ and
$\Xi(\Gamma_{\mathcal {B}}, \Delta)$ for  $r = 1$, $2$, $b-2$ and
$b-1$, respectively, where $\Delta$ is a self-paired $X$-orbit on
$Arc_3(\Gamma_{\mathcal {B}})$ and $\Lambda$ is a self-paired
$X$-orbit on $J(\Gamma_{\mathcal {B}})$. This  motivates us in this
section to consider the general case where $1 \leq r \leq b - 1$,
and introduce the stars and the double stars for a given graph. We
shall show that there is a close connection between $\Gamma$ and the
graph constructed from a certain set of double stars of
$\Sigma:=\Gamma_{\mathcal {B}}$.

\subsection{Stars of symmetric graphs}
Let $\Sigma$ be an $X$-symmetric graph with valency no less that
$2$. For $\tau \in V(\Sigma)$ and an $\mbk$-subset $S$ of
$\Sigma(\tau)$, we call $\mfs(\tau,S):=\{(\tau, \sigma) \in
Arc(\Sigma)\div \sigma \in S\}$ a {\em $\mbk$-star} of $\Sigma$ with
respect to $\tau$ and $S$. Set $\mcS t_{\tau}^{\mbk}(\Sigma) :=
\{\mfs(\tau,S)\div S \subseteq \Sigma(\tau), |S| = \mbk\}$ and $\mcS
t^{\mbk}(\Sigma) := {\cup_{\tau \in V(\Sigma)}}\mcS
t_{\tau}^{\mbk}(\Sigma).$ A star $\mfs:=\mfs(\tau,S)$ is said to be
{\em $X_{\mfs}$-symmetric} if $X_\tau\cap X_S$ acts transitively on
$S$. A nonempty subset $\mathcal{S}$ of $\mcS t^{\mbk}(\Sigma)$ is
said to be {\em $X$-symmetric} if $\mcS$ is $X$-transitive and
$\mfs$ is $X_{\mfs}$-symmetric for some $\mfs\in \mcS$.

Let $\mcS$ be an  $X$-symmetric subset of $\mcS t^{\mbk}(\Sigma)$.
For $\tau\in V(\Sigma)$, set  $\mcS_\tau=\{\mfs\in \mcS\div
\mfs=\mfs(\tau,S), S\subseteq\Sigma(\tau), |S|=\mbk\}$. Define an
incidence structure $\mathbb{D}(\tau) := (\Sigma(\tau), \mcS_{\tau},
\|)$ in which $\sigma \| \mfs (\tau ,S)$, for $\sigma \in
\Sigma(\tau)$, $\mfs (\tau ,S)\in S_{\tau}$, if and only if $\sigma
\in S$. A pair $(\sigma, \mfs)$ with $\sigma \| \mfs $ is said to be
a {\em flag} of $\mathbb{D}(\tau)$. Let $\mathbbm{r} := |\{\mfs
(\tau, S)\in \mcS_{\tau}\div \sigma \in S\}|, \mathbbm{b} :=
|\mcS_{\tau}|$ and $\mathbbm{v}:=val(\Sigma)$. Then it is easy to
see that $\mathbb{D}(\tau)$ is an $X_{\tau}$-flag-transitive
$1$-$(\mathbbm{v}, \mbk, \mathbbm{r})$ design with $\mathbbm{b}$
blocks. Moreover, $\mathbb{D}(\tau)$ is independent of the choice of
$\tau \in V(\Sigma)$ up to isomorphism.

The following Lemma~\ref{lem_star_design} says that, for $\tau\in
V(\Sigma)$, an arbitrary $X_{\tau}$-flag-transitive
$1$-$(\mathbbm{v}, \mbk, \mathbbm{r})$ design can be constructed as
above in some sense. Let
$\mathfrak{D}(\tau):=(\Sigma(\tau),\mathfrak{B}, \rm{I})$ be an
$X_{\tau}$-flag-transitive $1$-$(\mathbbm{v}, \mbk, \mathbbm{r})$
design. It may happen that distinct blocks $\mathfrak{b}_1$ and
$\mathfrak{b}_2$ of $\mathfrak{D}(\tau)$ have the same trace
$\{\sigma\div \sigma \rm{I}\mathfrak{b}_1\}=\{\sigma\div \sigma
\rm{I}\mathfrak{b}_2\}$. Since $\mathfrak{D}(\tau)$ is
flag-transitive, the number of blocks with the same trace is a
constant, say $m (\mathfrak{D}(\tau))$, called the multiplicity of
$\mathfrak{D}(\tau)$. Let $\mathfrak{D}'(\tau)$ be the design with
vertex set $\Sigma(\tau)$ and blocks being the traces of blocks of
$\mathfrak{D}(\tau)$. Then $\mathfrak{D}'(\tau)$ is an
$X_\tau$-flag-transitive $1$-$(\mathbbm{v},\mbk, \mathbbm{r}')$
design, where $\mathbbm{r}'=\frac{\mathbbm{r}}{m
(\mathfrak{D}(\tau))}$.
\begin{lem}\label{lem_star_design}
Let $\Sigma$ be an $X$-symmetric graph with valency $\mathbbm{v}
\geq 2$ and $\mathfrak{D}(\tau)$ be an $X_{\tau}$-flag-transitive
$1$-$(\mathbbm{v}, \mbk, \mathbbm{r})$ design with ${\mathbbm{b}}$
blocks, where $1 \leq \mbk \leq \mathbbm{v}-1$ and $\tau \in
V(\Sigma)$. Set $\mcS:=\{\mfs (\tau^x, S^x)\div x\in X, S\in
\mathfrak{D}'(\tau)\}$. Then $\mcS$ is $X$-symmetric, and
$\mathfrak{D}'(\tau)\cong \mathbb{D}(\tau)$ is an
$X_\tau$-flag-transitive $1$-$(\mathbbm{v},\mbk,
\frac{\mathbbm{r}}{m (\mathfrak{D}(\tau))})$ design with
$\frac{\mathbbm{b}}{m (\mathfrak{D}(\tau))}$ blocks.
\end{lem}

\subsection{Double stars}
Let $L$ and $R$ be $\mbk$-subsets of $\Sigma(\tau)$ and
$\Sigma(\sigma)$ respectively, set
$\mathfrak{l}=\mathfrak{s}(\tau,L)$ and
$\mathfrak{r}=\mathfrak{s}(\sigma,R)$, the pair $(\mfl,\mfr)$ is
called a {\em $\mbk$-double star} of $\Sigma$ if  $\sigma \in L$ and
$\tau \in R$. Denote by $D\mcS t^{\mbk}(\Sigma)$ the set of
$\mbk$-double stars of $\Sigma$. A nonempty subset $\Theta$ of
$D\mcS t^{\mbk}(\Sigma)$ is said to be {\em $X$-symmetric} if
$\mathcal{S}t(\Theta):=\{\mfl,\mfr\div (\mfl,\mfr) \in \Theta\}$ is
$X$-symmetric; and is self-paired if $(\mfl,\mfr) \in \Theta$ always
implies that $(\mfr,\mfl) \in \Theta$.

Here we give a straightforward lemma by ignoring the proof.
\begin{lem}\label{lem_symmetric_doublestar}
Let $\Sigma$ be an $X$-symmetric graph with valency $\mathbbm{v}
\geq 2$ and $\mbk$ an integer with $1 \leq \mbk \leq \mbv$.
\begin{enumerate}
 \item[{(a)}] If $\mcS$
is an $X$-symmetric orbit on $\mcS t^{\mbk}(\Sigma)$, then for
$\mfl=\mfs(\tau,L), \mfr=\mfs(\sigma,R) \in \mcS$ with $\sigma\in L$
and $\tau\in R$, $\Theta :=\{(\mfl^x, \mfr^x)\div x\in X\}$ is an
$X$-symmetric orbit on $D\mcS t^{\mbk}(\Sigma)$ and  $\mcS t(\Theta)
= \mcS$.

\item[{(b)}]Let $\Theta$ be an $X$-symmetric orbit  on $D\mcS
t^{\mbk}(\Sigma)$ and let $\tau, \sigma \in V(\Sigma)$. Then $(\tau,
\sigma) \in Arc(\Sigma)$ if and only if there exist
$\mfl:=\mfs(\tau,L), \mfr:=\mfs(\sigma, R) \in \mcS t(\Theta)$ such
that $(\mfl, \mfr) \in \Theta$.
\end{enumerate}
\end{lem}

The following example shows that an $X$-symmetric orbit $\Theta$ of
$\mbk$-double stars of an $X$-symmetric graph is not necessarily
self-paired.
\begin{ex}
{\rm Let $\Sigma$ be a cubic $(X,2)$-arc-regular graph with a
$3$-arc $(\tau_1, \tau, \sigma, \sigma_1)$  such that there is no
$x\in X$ maps this $3$-arc into $(\sigma_1, \sigma, \tau, \tau_1)$.
(See \cite[18c]{AGT}, for example.) Set $L= \{\tau_1, \sigma\}$, $R=
\{\sigma_1, \tau\}$, $\mfl := \mfs(\tau,L)$ and $\mfr= \mfs(\sigma,
R) $. Let $\Theta=\{(\mfl^x,\mfr^x)\div x\in X\}$. Then $\Theta$ is
an $X$-symmetric orbit on $D\mcS t^2(\Sigma)$. However, it is easily
shown that $\Theta$ is not self-paired.}
\end{ex}

\begin{con}\label{con_doublestar_graph}
Let $\Sigma$ be an $X$-symmetric graph with valency $\mathbbm{v}
\geq 2$ and $\Theta$ be a self-paired $X$-symmetric orbit on $D\mcS
t^{\mbk}(\Sigma)$ with $1 \leq \mbk \leq \mathbbm{v}-1$. Define a
graph $\Pi(\Sigma, \Theta)$, called the double star graph of
$\Sigma$ with respect to $\Theta$, with vertex set $St(\Theta)$ such
that two $\mbk$-stars $\mfl$ and $\mfr$ in $\mcS t(\Theta)$ are
adjacent if and only if $(\mfl, \mfr) \in \Theta$.
\end{con}

\begin{thm}\label{thm_doublestar_graph}
Let $\Sigma$, $\Theta$ and $\Gamma:=\Pi(\Sigma, \Theta)$ be as in
Construction~{\rm \ref{con_doublestar_graph}}. Set $\mcS= \mcS
t(\Theta)$ and $\mcB= \{\mcS_{\tau}\div \tau \in V(\Sigma)\}$, where
$\mcS_\tau=\{\mfs\in \mcS\div \mfs=\mfs(\tau,S),
S\subseteq\Sigma(\tau), |S|=\mbk\}$. Then $(\Gamma, X, \mathcal {B})
\in \mathcal {G}$ such that $\Gamma_{\mathcal {B}} \cong \Sigma$,
and for $B=\mcS_\tau\in \mcB$, $\mathcal {D}(B) \cong
\mathbb{D}^*(\tau)$, where $\mathbb{D}^*(\tau)$ is the dual design
of $\mathbb{D}(\tau)$.
\end{thm}
\proof It is easy to see that $\mathcal {B}$ is an $X$-invariant
partition of $ V(\Gamma)=\mcS$. For any $\mfs:=\mfs(\tau,S) \in
B:=\mcS_{\tau} \in \mathcal {B}$, as $1\le \mbk=|S| \le \mbv-1$,
take $\sigma \in S$ and $\delta \in \Sigma(\tau) \setminus S$. Since
$\Sigma$ is $X$-symmetric, there exists $g \in X_{\tau}$ such that
$\delta = \sigma^g$. Let $\mfl= \mfs^g$. Then $\mfs\ne \mfl\in
\mcS_{\tau}$, thus $v = |\mcS_{\tau}| \geq 2$, and hence $\mathcal
{B}$ is a nontrivial $X$-invariant partition of $ V(\Gamma)$. By
Lemma \ref{lem_symmetric_doublestar}, there exists $\mfr \in
\mcS_{\delta}$ such that $(\mfl, \mfr) \in \Theta$, hence
$C:=\mcS_{\delta} \in \Gamma_{\mathcal {B}}(B)$. If there exists
$\mfr' \in \mcS_{\delta}$ such that $(\mfs, \mfr') \in \Theta$, then
$\delta \in S$, a contradiction. Thus $\mfs \not\in B\cap \Gamma(C)$
and $k=|B\cap \Gamma(C)| \leq v-1$. Hence $\Gamma$ is not a
multicover of $\Gamma_{\mathcal {B}}$. It is easily shown by using
Lemma \ref{lem_symmetric_doublestar} that $\Gamma$ is $X$-symmetric,
and $V(\Sigma) \rightarrow V(\Gamma_{\mathcal {B}}),\, \tau \mapsto
\mcS_{\tau}$ is an isomorphism from $\Sigma$ to $\Gamma_{\mathcal
{B}}$.

For $\tau \in V(\Sigma)$ and $B=\mcS_\tau$, define a map $\pi: B\cup
\Gamma_\mcB(B)\rightarrow \mcS_\tau\cup \Sigma(\tau);\, \mfs(\tau,
S)\mapsto \mfs(\tau, S),\, C\mapsto \sigma$ for $\mfs(\tau, S)\in
B=\mcS_\tau$ and $C=\mcS_\sigma\in \Gamma_\mcB(B)$. Assume
$C=\mcS_\sigma\in \Gamma_\mcB(B)$. Then by the definition of
$\Gamma_\mcB$ and the construction of $\Gamma$ there exist
$\mfl=\mfs(\tau,L)\in B$ and $\mfr=\mfs(\sigma,R)\in C$ such that
$(\mfl,\mfr)\in \Theta$. In particular, $\sigma\in L\subseteq
\Sigma(\tau)$. Thus $\pi$ is well-defined. Moreover $\pi$ is a
bijection. By the definition of $\mathcal{D}(B)$, for
$\mfs=\mfs(\tau,S)\in B$ and $ C=\mcS_\sigma\in \Gamma_\mcB(B)$, we
know that $\mfs{\rm |} C$ if and only if there is some
$\mathfrak{t}=\mfs(\sigma,T)\in C$ such that $(\mfs,\mathfrak{t})\in
\Theta$, that is, $\tau\in T$ and $\sigma \in S$; it follows  that
$\sigma\| \mfs$. Now assume that $\sigma'\in \Sigma(\tau)$ and
$\mfs'=\mfs(\tau,S')$ with $\sigma'\| \mfs'$. Then $\sigma'\in S'$.
Take some $\mathfrak{t}'=\mfs(\tau', T')$ such that
$(\mfs',\mathfrak{t}')\in \Theta$. Then $\tau'\in S'$. Since $\mfs'$
is $X_{\mfs'}$-symmetric, there is some $x\in X_\tau\cap X_{S'}$
with $\tau'^x=\sigma'$. Thus $\mfs'^x=\mfs'$,
$\mathfrak{t}'^x=\mfs(\sigma', T'^x)\in \mcS_{\sigma'}$ and
$(\mfs',\mathfrak{t}'^x)=(\mfs',\mathfrak{t}')^x\in \Theta$. Hence
$\mfs'{\rm |} \mcS_{\sigma'}$. The above argument says that $\pi$ is
an isomorphism from $\mathcal {D}(B)$ to $\mathbb{D}^*(\tau)$. So
$\mathcal {D}(B) \cong \mathbb{D}^*(\tau)$. \qed

Here we give the following sufficient condition which is useful in
determining whether or not a double star graph exists.
\begin{thm}\label{thm_design_odd}
Let $\Sigma$ be an $X$-symmetric graph with valency $\mathbbm{v}
\geq 2$ and let $\tau \in V(\Sigma)$. If there exists some
$X_{\tau}$-flag-transitive $1$-$(\mathbbm{v}, \mbk, \mathbbm{r})$
design $\mfD(\tau)$ on $\Sigma(\tau)$ for $1\le \mbk\le \mbv-1$ such
that $\frac{\mathbbm{r}}{m (\mathfrak{D}(\tau))}$ is odd, then there
exists a self-paired $X$-symmetric orbit $\Theta$ on $D\mcS
t^{\mbk}(\Sigma)$.
\end{thm}
\proof By Lemma \ref{lem_star_design}, setting $\mcS=\{\mfs (\tau^x,
S^x)\div x\in X, S\in \mathfrak{D}'(\tau)\}$, we know that
$\mathfrak{D}'(\tau)\cong\mathbb{D}(\tau)$  is an
$X_\tau$-flag-transitive $1$-$(\mathbbm{v},\mbk,
\frac{\mathbbm{r}}{m (\mathfrak{D}(\tau))})$ design with
$\frac{\mathbbm{b}}{m (\mathfrak{D}(\tau))}$ blocks, and $\mcS$ is
$X$-symmetric. Let $(\tau, \sigma)\in Arc(\Sigma)$. Then, since
$\Sigma$ is $X$-symmetric, $(\tau, \sigma)^y = (\sigma, \tau)$ for
some $y \in X$. Set $\mcS_{(\tau,\sigma)}=\{\mfs(\tau, S)\in
\mcS_\tau\div\sigma\in S\}$. Then $\frac{\mathbbm{r}}{m
(\mathfrak{D}(\tau))}=|\mcS_{(\tau,\sigma)}|$ is odd,
$\mcS_{(\tau,\sigma)}^y=\mcS_{(\sigma,\tau)}$ and
$\mcS_{(\tau,\sigma)}^{y^2}=\mcS_{(\tau,\sigma)}$. Let $\mathcal{O}$
be a $\langle y^2 \rangle$-orbit  on $S_{(\tau, \sigma)}$ with odd
length $l$. Then for $\mfl\in \mathcal{O}$, the stabilizer of $\mfl$
in  $\langle y^2\rangle $ is  $\langle y^{2l}\rangle$. Let $z=y^l$
and $\mfr=\mfl^z$. Then $(\mfl,\mfr)^z=(\mfr,\mfl)$, and hence
$\Theta := \{(\mfl^x, \mfr^x)\div x\in X\}$ is a self-paired
$X$-symmetric orbit on $D\mcS t^\mbk(\Sigma)$ with  $\mcS t(\Theta)=
\mcS$. \qed

The following Theorem~\ref{thm_quotient_star_design} says that, for
any $X$-symmetric graph $\Gamma$ with an nontrivial $X$-invariant
partition, $\Gamma$ or  a quotient of $\Gamma$ can be constructed as
in Construction~\ref{con_doublestar_graph}.

Let $(\Gamma, X, \mathcal {B}) \in \mathcal {G}$. For $B \in
\mathcal {B}$ and $\mathfrak{v}\in B$, define
$B_\mfv=B\cap(\cap_{C\in \Gamma_\mathcal{B}(\mfv)}\Gamma(C))$. Then
$|B_\mfv|$, denoted by $m^*(\Gamma,\mcB)$ is independent of the
choices of $B$ and $\mfv$. Noting that $\Gamma$ is not a multicover
of $\Gamma_\mcB$, we have $m^*(\Gamma,\mcB)\le k:=|B\cap \Gamma(C)|$
for $C\in \Gamma_\mcB(B)$. In fact, $ m^*(\Gamma,\mcB)$ is  the
multiplicity of the dual design $\mathcal{D}^*(B)$ of
$\mathcal{D}(B)$. Set $\underline{\mcB}=\{B_\mfv\div B\in \mcB,
\mfv\in B\}$. Then $\underline{\mcB}$ is an $X$-invariant partition
of $V(\Gamma)$. For $B\in \mcB$, we set $\bar{B}=\{B_\mfv\div
\mfv\in B\}$. Then $\Gamma_{\underline{\mcB}}$ is an $X$-symmetric
graph with an $X$-invariant partition $\bar{\mcB}:=\{\bar{B}\div
B\in \mcB\}$ such that
$(\Gamma_{\underline{\mcB}})_{\bar{\mcB}}\cong \Gamma_\mcB$.
Moreover, $m^*(\Gamma_{\underline{\mcB}},\bar{\mcB})=1$.

\begin{thm}\label{thm_quotient_star_design}
Let $(\Gamma, X, \mathcal {B}) \in \mathcal {G}$. Set $\mcS =
\{\mfs(B, \Gamma_\mcB(\mfv))\div B\in \mcB, \mfv\in B\}$. Then
$\mcS$ is an $X$-symmetric orbit on $\mcS t^{r}(\Gamma_{\mathcal
{B}})$, where $r=|\Gamma_\mcB(\mfv)|$ is a constant. Let $\Theta=
\{(\mfl, \mfr)\div \mfl=\mfs(B, \Gamma_\mcB(\mfv)), \mfr=\mfs(C,
\Gamma_\mcB(\mfu)), \mfv\in B\in \mcB, \mfu\in C\in \mcB, (\mfv,
\mfu) \in Arc(\Gamma)\}$. Then $\Theta$ is a self-paired
$X$-symmetric orbit on $D\mcS t^{r}(\Gamma_{\mathcal {B}})$ with
$\mcS t(\Theta) = \mcS$ and $\Gamma_{\underline{\mcB}} \cong
\Pi(\Gamma_{\mathcal {B}}, \Theta)$, and $X$ acts faithfully on
$\mcB$ if and only $X$ acts faithfully on $\underline{\mcB}$.
\end{thm}
\proof It is easily shown that $\Theta$ is a self-paired
$X$-symmetric orbit on $D\mcS t^{r}(\Gamma_{\mathcal {B}})$ with
$\mcS t(\Theta) =\mcS$. Assume $ m^*(\Gamma,\mcB)=1$. Then, for two
distinct vertices $\mfv\in B\in \mcB$ and $\mfu\in C\in \mcB$ of
$\Gamma$, $B_\mfv=\{\mfv\}$ and $C_\mfu=\{\mfu\}$, it implies
$\Gamma_\mcB(\mfv)\ne\Gamma_\mcB(\mfu)$, and hence  $\mfs(B,
\Gamma_\mcB(\mfv))\ne\mfs(C, \Gamma_\mcB(\mfu))$. Thus
$V(\Gamma)\rightarrow V(\Pi(\Gamma_{\mathcal {B}})), \, \mfv\mapsto
\mfs(B, \Gamma_\mcB(\mfv))$ is a bijection. Further, it is easy to
see this bijection is in fact an isomorphism between $\Gamma$ and
$\Pi(\Gamma_{\mathcal {B}},\Theta)$.

Now assume $ m^*(\Gamma,\mcB)>1$. Recall that $ m^*(\Gamma,\mcB)\le
k:=|B\cap \Gamma(C)|$ for $C\in \Gamma_\mcB(B)$. Then
$\underline{\mcB}$ is a proper refinement of $\mcB$. Consider
$\Gamma_{\underline{\mcB}}$ with $X$-invariant partition
$\bar{\mcB}$. Then $ m^*(\Gamma_{\underline{\mcB}},\bar{\mcB})=1$.
Then a similar argument as above leads to $\Gamma_{\underline{\mcB}}
\cong \Pi(\Sigma, \bar{\Theta})$, where
$\Sigma=(\Gamma_{\underline{\mcB}})_{\bar{\mcB}}$ and
$\bar{\Theta}=\{(\bar{\mfl}, \bar{\mfr})\div
\bar{\mfl}=\mfs(\bar{B}, \Sigma(B_\mfv)), \bar{\mfr}=\mfs(\bar{C},
\Sigma(C_\mfu)), B_\mfv\in \bar{B}\in \bar{\mcB}, C_\mfu\in
\bar{C}\in \bar{\mcB}, (B_\mfv, C_\mfu) \in
Arc(\Gamma_{\underline{\mcB}})\}$. Noting that $B_\mfv=B_{\mfv'}$
for any $\mfv'\in B_\mfv$, it follows that $\mfs(\bar{B},
\Sigma(B_\mfv))\mapsto \mfs(B, \Gamma_\mcB(\mfv))$ gives a bijection
between $V(\Pi(\Sigma, \bar{\Theta}))$ and $V(\Pi(\Gamma_{\mathcal
{B}},\Theta))$, which is in fact an isomorphism between $\Pi(\Sigma,
\bar{\Theta})$ and $\Pi(\Gamma_{\mathcal {B}},\Theta)$. Hence
$\Gamma_{\underline{\mcB}} \cong \Pi(\Gamma_{\mathcal {B}},
\Theta)$.

Let $K$ and $H$ be the kernels of $X$ acting on $\mcB$ and on
$\underline{\mcB}$ respectively. Noting that $\underline{\mcB}$ is a
refinement of $\mcB$, we have $H\le K$. Let $x\in K$ and $B_\mfv\in
\bar{B}\in \underline{\mcB}$. Since $
m^*(\Gamma_{\underline{\mcB}},\bar{\mcB})=1$, we have
$\{B_\mfv\}=\bar{B}\cap(\cap_{\bar{C}\in
(\Gamma_{\underline{\mcB}})_{\bar{\mcB}}(B_\mfv)}\Gamma_{\underline{\mcB}}(\bar{C}))
=\bar{B}\cap(\cap_{C\in
\Gamma_\mcB(\mfv)}\Gamma_{\underline{\mcB}}(\bar{C}))$, yielding
$B_\mfv^x=B_\mfv$. The above argument implies $x\in H$. Hence $K\le
H$, and so $H=K$. Therefore,  $X$ acts faithfully on $\mcB$ (that
is, $K=1$) if and only $X$ acts faithfully on $\underline{\mcB}$
(that is, $H=1$). \qed

Finally, we list a simple fact which will be used in the following
sections.
\begin{thm}\label{thm_B_neighbour}
Let $(\Gamma, X, \mathcal {B}) \in \mathcal {G}$ and $B\in \mcB$. If
$ m^*(\Gamma,\mcB)=1$ and $m (\mathcal{\mathcal{D}}(B))=1$, then
$X_B^B \cong X_B^{\Gamma_{\mathcal {B}}(B)}$.
\end{thm}
\proof  If $x\in X$ fixes $B$ set-wise, then it also fixes the
neighborhood $\Gamma_\mcB(B)$ of $B$ in $\Gamma_\mcB$. Now consider
the action of $X_B$ on $\Gamma_\mcB(B)$, and let $K$ be the kernel
of this action. For any $\mfv\in B$, since $ m^*(\Gamma,\mcB)=1$, we
have $\{\mfv\}=B\cap(\cap_{C\in \Gamma_\mcB(\mfv)}\Gamma(C))$. It
follows that $K$ fixes $\mfv$. Thus $K\le X_{(B)}$. On the other
hand, $x$ fixes $B\cap\Gamma(C)$ point-wise for any $x\in X_{(B)}$
and any $C\in \Gamma_\mcB(B)$, in particular,
$B\cap\Gamma(C^x)=(B\cap\Gamma(C))^x=B\cap\Gamma(C)$. It follows
from $m (\mathcal{\mathcal{D}}(B))=1$ that $C=C^x$. Therefore, $x\in
K$. Thus $X_{(B)}\le K$, and so $X_{(B)}=K$. Then $X_B^B \cong
X_B/X_{(B)}=X_B/K \cong X_B^{\Gamma_{\mathcal {B}}(B)}$. \qed

\section{The main result}
We state the main result of this paper in this section and prove it
in the next four sections.

To state the result we need the following concept. A {\em near
$n$-gonal graph}~\cite{NG} is a connected graph $\Sigma$ of girth at
least $4$ together with a set $\mathcal{E}$ of $n$-cycles of
$\Sigma$ such that each $2$-arc of $\Sigma$  is contained in a
unique member of $\mathcal{E}$.

Let $(\Gamma, X, \mathcal {B}) \in \mathcal {G}$. For a subgraph
$\Delta$ of $\Gamma$, denote by $X_{[\Delta]}$ the subgroup of $X$
which preserves the adjacency of $\Delta$, and set
$X_{[\Delta]}^{[\Delta]}=X_{[\Delta]}/X_{(V(\Delta))}$. Recall that,
for $\mfv\in B\in \mcB$ and $C\in \Gamma_\mcB(B)$, the parameters $v
:= |B|$, $k := |\Gamma(C) \cap B|$, $r := |\Gamma_{\mathcal
{B}}(\mfv)|$ and $b := val(\Gamma_{\mathcal {B}})$ are independent
of the choices of $B$ and $\mfv$, and $\mathcal {D}(B)$ is an
$X_B$-flag-transitive $1$-$(v, k, r)$ design with $b$ blocks. Now we
are ready to state the main result of this paper.
\begin{thm}\label{thm_main}
Let $(\Gamma, X, \mathcal {B}) \in \hat{\mathcal {G}}$ and $B \in
\mathcal {B}$. Let $e= |E(\Gamma_{\mathcal {B}})|$, $\mu
=|V(\Gamma_{\mathcal {B}})|$. Suppose that $k = 3$. Then
$\Gamma_{\mathcal {B}}$ is $(X,2)$-arc-transitive if and only if one
of the following four cases occurs.
\begin{enumerate}
\item[{(a)}] $(v, b, r) = (4, 4, 3)$, $X_B^B \cong A_4$ or
$S_4$;
\item[{(b)}] $(v, b, r) = (6, 4, 2)$,
$X_B^B \cong A_4$ or $S_4$;
\item[{(c)}] $(v, b, r) = (7, 7, 3)$,
 $X_B^B \cong PSL(3,2)$;
\item[{(d)}] $v = 3b\geq 6$, $r = 1$ and
$X_B$ acts $2$-transitively on the blocks of $\mathcal {D}(B)$.
\end{enumerate}
Furthermore, if case (a) occurs, then $\Gamma \cong
\Xi(\Gamma_{\mathcal {B}}, \Delta)$ for some self-paired $X$-orbit
$\Delta$ on $Arc_3(\Gamma_{\mathcal {B}})$, $X$ acts faithfully on
$\mcB$, and any connected tetravalent $(X,2)$-arc-transitive graph
can occur as $\Gamma_{\mathcal {B}}$; moreover, one of the following
three statements holds.
\begin{enumerate}
\item[(a.1)]
$\Gamma[B, C] \cong 3K_2$,  $val(\Gamma) = 3$, there exists an
$X$-orbit $\mathcal {E}$ of $n$-cycles of $\Gamma_\mcB$ with
$|\mathcal {E}| = m$ such that $\Delta=\cup_{C\in \mcE}Arc_3(C)$,
$X_{[C]}^{[C]} \cong D_{2n}$ for each $C\in \mathcal {E}$, where $m
\geq 6$ and $n \geq girth(\Gamma_\mcB)$ with $mn = 3e = 6\mu$.
Moreover, either $\Gamma_\mcB \cong K_{5}$ or $\Gamma_\mcB$ is a
near $n$-gonal graph with respect to $\mathcal {E}$;  either $X_B
\cong A_4$,  $\Gamma$ is  $(X,1)$-arc-regular and $\Gamma_\mcB$ is
$(X,2)$-arc-regular; or $X_B  \cong S_4$ and $\Gamma$ is
$(X,2)$-arc-regular.
\item[(a.2)]
$\Gamma[B, C] \cong K_{3,3} - 3K_2$, $val(\Gamma) = 6$, $X_B\cong
S_4$, and $\Gamma$ is connected and $(X,1)$-arc-regular. Further,
$\Delta' := Arc_3(\Gamma_{\mathcal {B}})\setminus \Delta$ is a
self-paired $X$-orbit on $Arc_3(\Gamma_{\mathcal {B}})$, and there
exists an $X$-orbit $\mathcal {E}$ of $n$-cycles of $\Sigma$ with
$|\mathcal {E}| = m$, such that $\Delta' = \cup_{C \in \mathcal
{E}}Arc_3(C)$, $X_{[C]}^{[C]} \cong D_{2n}$ for each $C \in \mathcal
{E}$, where $m \geq 6$ and $n \geq girth(\Sigma)$ with $mn = 3e =
6\mu$. Moreover, either $\Gamma_{\mathcal {B}} \cong K_5$ or
$\Gamma_{\mathcal {B}}$ is a near $n$-gonal graph.
\item[(a.3)]
$\Gamma[B, C] \cong K_{3,3}$, $val(\Gamma) = 9$, $\Gamma$ is
connected and $(X,1)$-transitive, and $\Gamma_\mcB$ is $(X,3)$-arc
transitive.
\end{enumerate}
If case (b) holds, then $\Gamma  \cong \gimel(\Gamma_{\mathcal {B}},
\Delta) \cong \Psi{(\Gamma_{\mathcal {B}}}, \Lambda)$ for some
self-paired $X$-orbit $\Delta$ on $Arc_3(\Gamma_{\mathcal {B}})$ and
some self-paired $X$-orbit $\Lambda$ on $J(\Gamma_\mcB)$, $X$ acts
faithfully on ${\mathcal {B}}$, and any connected tetravalent
$(X,2)$-arc-transitive graph can occur as $\Gamma_{\mathcal {B}}$;
moreover, one of the following three cases occurs.
\begin{enumerate}
\item[(b.1)]
$\Gamma[B,C] \cong 3K_2$, $\Gamma\cong mC_n$,
 and there exists an $X$-orbit $\mathcal {E}$ of
$n$-cycles of $\Gamma_\mcB$ with $|\mathcal {E}| = m$, such that
$\Delta=\cup_{C\in \mcE}Arc_3(C)$, $X_{[C]}^{[C]} \cong D_{2n}$ for
each $C\in \mathcal {E}$, where $m \geq 6$ and $n \geq
girth(\Gamma_\mcB)$ with $mn = 3e = 6\mu$. Moreover, either
$\Gamma_\mcB \cong K_{5}$ or $\Gamma_\mcB$ is a near $n$-gonal graph
with respect to $\mathcal {E}$; either $X_B\cong A_4$, $\Gamma$ is
$(X,1)$-arc-regular and $\Gamma_\mcB$ is $(X,2)$-arc-regular, or
$X_B\cong S_4$ and $\Gamma$ is not $(X,1)$-arc-regular.
\item[(b.2)]
$\Gamma[B, C] \cong K_{3,3} - 3K_2$, $val(\Gamma) = 4$, $X_B\cong
S_4$, $\Gamma$ is connected and $(X,1)$-arc-regular. Further,
$\Delta' := Arc_3(\Gamma_{\mathcal {B}})\setminus \Delta$ is a
self-paired $X$-orbit on $Arc_3(\Gamma_{\mathcal {B}})$, and there
exists an $X$-orbit $\mathcal {E}$ of $n$-cycles of $\Sigma$ with
$|\mathcal {E}| = m$, such that $\Delta' = \cup_{C \in \mathcal
{E}}Arc_3(C)$, $X_{[C]}^{[C]} \cong D_{2n}$ for each $C \in \mathcal
{E}$, where $m \geq 6$ and $n \geq girth(\Sigma)$ with $mn = 3e =
6\mu$. Moreover, either $\Gamma_{\mathcal {B}} \cong K_5$ or
$\Gamma_{\mathcal {B}}$ is a near $n$-gonal graph.
\item[(b.3)]
$\Gamma[B, C] \cong K_{3,3}$, $val(\Gamma) = 6$,  $\Gamma$ is
connected and $(X,1)$-transitive, and $\Gamma_\mcB$ is $(X,3)$-arc
transitive.
\end{enumerate}
If  case $(c)$ holds, then $\Gamma \cong \Pi(\Gamma_{\mathcal {B}},
\Theta)$ for some self-paired $X$-symmetric orbit $\Theta$ on $D\mcS
t^{3}(\Gamma_{\mathcal {B}})$, $X$ acts faithfully on $\mcB$, one
connected heptavalent $X$-symmetric graph $\Sigma$ can occur as
$\Gamma_{\mathcal {B}}$, if and only if $X_\tau^{\Sigma(\tau)}\cong
PSL(3,2)$ for $\tau \in V(\Sigma)$; further, one of the following
three cases occurs.
\begin{enumerate}
\item[(c.1)] $\Gamma[B,C] \cong 3K_2$, $val(\Gamma) = 3$, $\Gamma$ is
$(X,2)$-arc-transitive but not $(X,2)$-arc-regular.

\item[(c.2)] $\Gamma[B,C] \cong K_{3,3}-3K_2 $, $val(\Gamma) = 6$,
$\Gamma$ is connected and is $(X,1)$-transitive.

\item[(c.3)] $\Gamma[B,C] \cong K_{3,3}$, $val(\Gamma) = 9$,
$\Gamma$ is connected and is $(X,1)$-transitive.
\end{enumerate}
If case (d) occurs, then one of the following three cases occurs.
\begin{enumerate}
\item[(d.1)] $\Gamma[B,C]\cong 3K_2$, $\Gamma \cong 3eK_2$.
\item[(d.2)] $\Gamma[B,C]\cong K_{3,3}-3K_2$, $\Gamma \cong
e(K_{3,3}-3K_2)$.
\item[(d.3)] $\Gamma[B,C]\cong K_{3,3}$, $\Gamma \cong eK_{3,3}$.
\end{enumerate}
\end{thm}

\section{Self-paired orbits of $3$-arcs}
We begin this section by showing that there always exists a
self-paired $X$-orbit  of $3$-arcs for any symmetric graph of even
valency.
\begin{thm}\label{thm_evenvalency_selfpaired}
Any $X$-symmetric graph $\Sigma$ of even valency $\mathbbm{v} \geq
2$ contains a self-paired $X$-orbit $\Delta$ on $Arc_3(\Sigma)$.
\end{thm}
\proof For any $(\tau, \sigma) \in Arc(\Sigma)$, as $\Sigma$ is
$X$-symmetric, there exists $y \in X$ such that $(\tau, \sigma)^y =
(\sigma, \tau)$, and so $(\Sigma(\tau) \setminus \{\sigma\})^y =
\Sigma(\sigma) \setminus \{\tau\}$, $(\Sigma(\tau) \setminus
\{\sigma\})^{y^2} = \Sigma(\tau) \setminus \{\sigma\}$. Since
$|\Sigma(\tau) \setminus \{\sigma\}| = \mathbbm{v}-1$ is odd, there
must be some $\langle y^2 \rangle$-orbit $\mathcal{O}$ on
$\Sigma(\tau) \setminus \{\sigma\}$ with odd length $l$. For
$\tau_1\in \mathcal{O}$, the stabilizer of $\tau_1$ in $\langle
y^2\rangle$ is $\langle y^{2l}\rangle$. Let $z = y^{l}$, $\sigma_1 =
\tau_1^z$ and $\bbalpha = (\tau_1, \tau, \sigma, \sigma_1)$. Since
$l$ is odd, $(\tau, \sigma)^z =(\tau, \sigma)^{y^l}= (\sigma,
\tau)$. Then $\bbalpha \in Arc_3(\Sigma)$ and   $\bbalpha^{z} =
\bbalpha^{-1}$. Thus $\Delta = \{(\tau_1^x, \tau^x, \sigma^x,
\sigma_1^x)\div x\in X\}$  is a self-paired $X$-orbit on
$Arc_3(\Sigma)$. \qed

Let $\Sigma$ be an $X$-symmetric graph with valency $\mathbbm{v}
\geq 2$ and $\Delta$ be an $X$-orbit on $Arc_3(\Sigma)$. For any
$\bbalpha := (\tau_1, \tau, \sigma, \sigma_1) \in \Delta$, consider
the action of $X_{(\tau_1, \tau, \sigma)}$  on $\Sigma(\sigma)
\setminus \{\tau\}$, and denote by $O_1, O_2,\ldots, O_t$ the orbits
of this action. Without loss of generality, assume $\sigma_1 \in
O_1$ and $|O_2| \leq |O_3| \leq \ldots \leq |O_t|$. Since $\Delta$
is an $X$-orbit of $3$-arcs, all $\ell_i(\Delta) := |O_i| \geq 1$
are independent of the choice of $\bbalpha \in \Delta$. Set
$\bfl(\Delta) = (\ell_1(\Delta), \ldots, \ell_t(\Delta))$.

\begin{thm}\label{thm_nearngonal-graph}
Let $\Sigma$ be a connected $(X,2)$-arc-transitive graph with
valency $\mathbbm{v} \geq 3$ and $\Delta$ be a self-paired $X$-orbit
on $Arc_3(\Sigma)$ such that $\ell_1(\Delta) = 1$. If  $X$ is
faithful on $V(\Sigma)$, then $X_\tau$ is faithful on $\Sigma(\tau)$
for $\tau\in V(\Sigma)$. Set $\mu = |V(\Sigma)|$ and $e =
|E(\Sigma)|$. Then $\gimel(\Sigma, \Delta) \cong mC_n$ such that
\begin{enumerate}
\item[(1)] $m \geq \mathbbm{v}(\mathbbm{v}-1)/2$, $n \geq girth(\Sigma) \geq 3$
and $mn = {\mu}{\mathbbm{v}}(\mathbbm{v}-1)/2 = e(\mathbbm{v} - 1)$;
\item[(2)] there exists an $X$-orbit $\mathcal {E}$ of
$n$-cycles  of $\Sigma$ with  $\Delta=\cup_{C\in \mcE}Arc_3(C)$ and
$|\mathcal {E}| = m$;
\item[(3)] $X_{[C]}^{[C]} \cong D_{2n}$ for  $C \in \mathcal
{E}$, where $D_{2n}$ is the dihedral group of order $2n$;
\item[(4)] every  $2$-path of $\Sigma$ is contained in a unique
member of $\mathcal {E}$, and either $\Sigma \cong
K_{\mathbbm{v}+1}$ {\rm (}the complete graph on $\mathbbm{v}+1$
vertices{\rm )}, or  $n\ge girth(\Sigma)\ge 4$ and $\Sigma$ is a
near $n$-gonal graph with respect to $\mathcal {E}$.
\end{enumerate}
\end{thm}
\proof  Since $\Sigma$ is $(X,2)$-arc-transitive, every $2$-arc of
$\Sigma$ lies in a member of $\Delta$. Let $(\tau,\sigma)$ be an
arbitrary arc of $\Sigma$. Since $\ell_1(\Delta)=1$ and $\Delta$ is
a self-paired $X$-orbit, we conclude that, for any $\tau_1\in
\Sigma(\tau)\setminus\{\sigma\}$ there is a unique $\sigma_1\in
\Sigma(\sigma)\setminus\{\tau\}$ such that
$(\tau_1,\tau,\sigma,\sigma_1)\in \Delta$,
$X_{(\tau_1,\tau,\sigma)}=X_{(\tau,\sigma,\sigma_1)}$, and that
$(\tau_1',\tau,\sigma,\sigma_1)\in \Delta$ implies $\tau'_1=\tau$.
Then $(X_\tau)_{(\Sigma(\tau))}=\cap_{\tau_1\in
\Sigma(\tau)\setminus\{\sigma\}}X_{(\tau_1,\tau,\sigma)}=\cap_{\sigma_1\in
\Sigma(\sigma)\setminus\{\tau\}}X_{(\tau,\sigma,\sigma_1)}=(X_\sigma)_{(\Sigma(\sigma))}$.
It follows from the connectedness of $\Sigma$ that
$(X_\tau)_{(\Sigma(\tau))}$ fixes every vertex of $\Sigma$. Thus, if
$X$ is faithful on $V(\Sigma)$, then $(X_\tau)_{(\Sigma(\tau))}=1$
and $X_\tau$ is faithful on $\Sigma(\tau)$.

Let $\gimel= \gimel(\Sigma, \Delta)$.  By
Proposition~\ref{prop_U-graph},   $\gimel$ is $X$-symmetric and
admits an $X$-invariant partition $\mathcal {P} := \{P_{\sigma} \div
\sigma \in V(\Sigma)\}$ such that $\Sigma\cong \gimel_\mathcal{P}$,
where $P_\sigma$ is the set of $2$-paths of $\Sigma$ with middle
vertex $\sigma$. It follows from \cite{FSGTTQ(II)} that
$r:=|\gimel_\mathcal{P}(\mfv)|=2$ and $\lambda :=|P_\delta\cap
\gimel(P_\tau) \cap  \gimel(P_\sigma)|=1$ for any vertex $\mfv$ (a
$2$-path of $\Sigma$) in $V(\gimel)$ and $P_\delta$ with $\mfv\in
P_\delta$ and $\gimel_\mathcal{P}(\mfv)=\{P_\tau,P_\sigma\}$. Since
$\ell_1(\Delta) = 1$ and $\Delta$ is self-paired, for any $2$-path
$[\tau_1, \tau, \sigma]$ of $\Sigma$, there exist exactly two
$2$-paths $[\tau, \sigma, \sigma_1]$ and $[\tau_2,\tau_1, \tau]$
such that $(\tau_1, \tau, \sigma, \sigma_1) \in \Delta$ and
$(\tau_2, \tau_1, \tau, \sigma) \in \Delta$. It follows that
$\gimel$ is of valency two, and that $\gimel$ is a disjoint union of
cycles. Assume $\gimel\cong mC_n$. Then $mn$ is the number of
$2$-paths of $\Sigma$, and hence
$mn={\mu}{\mathbbm{v}}(\mathbbm{v}-1)/2 = e(\mathbbm{v} - 1)$.
Noting that $\gimel$ is of valency $2$ and every $P_{\sigma}$ is an
independent set of $V(\gimel)$, it follows that different vertices
in $P_{\sigma}$ appear in different $n$-cycles of $\gimel$. Thus $m
\geq |P_{\sigma}|=\mathbbm{v}(\mathbbm{v}-1)/2$.

Let $\mathfrak{C}=[\mfv_1, \mfv_2, \ldots, \mfv_n, \mfv_1]$ be an
arbitrary $n$-cycle of $\gimel$, where $\mfv_i=[\tau_{i},
\sigma_{i}, \delta_{i}]$ are $n$ distinct $2$-paths of $\Gamma$ with
middle vertices $\sigma_i$, respectively. Without loss of
generality, we assume $\delta_i=\sigma_{i+1}=\tau_{i+2}$ for $1\le
i\le n$, where subscripts are reduced modulo $n$. Since $\mfv_i$ is
a $2$-path, $\sigma_i\ne\delta_i$, hence $\sigma_i\ne \sigma_{i+1}$.
Then $(\sigma_i,\sigma_{i+1})\in Arc(\Sigma)$. Since
$\{\mfv_i,\mfv_{i+1}\}$ is an edge of $\gimel$, we have
$(\sigma_{i-1}, \sigma_i, \sigma_{i+1},  \sigma_{i+2}
)=(\tau_i,\sigma_i,\delta_i,\delta_{i+1})\in \Delta$.

Now we shall show $C=[\sigma_1,\sigma_2, \ldots, \sigma_n,
\sigma_1]$ is an $n$-cycle of $\Sigma$. In particular, $n\ge
girth(\Sigma)\ge 3$. Note that $\mfC$ is a  component of $\gimel$.
Then $\mfC$ is $X_{[\mfC]}$-symmetric; in particular,
$X_{[\mfC]}^{[\mfC]}\cong D_{2n}$, the dihedral group of order $2n$.
Thus there exist $x,y\in X_{[\mfC]}$ such that $\mfv_i^x=\mfv_{i+1}$
and $\mfv_i^y=\mfv_{n-i+1}$, hence $\sigma_i^x=\sigma_{i+1}$ and
$\sigma_i^y=\sigma_{n-i+1}$ for $1\le i\le n$ with subscripts modulo
$n$. Assume that $\sigma_i=\sigma_j$ for some $i$ and $j$. Then
$\sigma_{i+1}=\sigma_i^x=\sigma_j^x=\sigma_{j+1}$ and
$\sigma_{i+2}=\sigma_{i+1}^x=\sigma_{j+1}^x=\sigma_{j+2}$. Thus
$P_{\sigma_i}=P_{\sigma_j}$, $P_{\sigma_{i+1}}=P_{\sigma_{j+1}}$ and
$P_{\sigma_{i+2}}=P_{\sigma_{j+2}}$. It yields
$(\mfv_i,\mfv_{i+1}),\,(\mfv_j,\mfv_{j+1})\in
Arc(\gimel[P_{\sigma_i},P_{\sigma_{i+1}}])$ and
$(\mfv_{i+1},\mfv_{i+2}),\,(\mfv_{j+1},\mfv_{j+2})\in
Arc(\gimel[P_{\sigma_{i+1}},P_{\sigma_{i+2}}])$. It follows that
$\mfv_{i+1},\, \mfv_{j+1}\in P_{\sigma_{i+1}}\cap \gimel(
P_{\sigma_i})\cap\gimel( P_{\sigma_{i+2}})$. Since $1=\lambda
=|P_{\sigma_{i+1}}\cap \gimel( P_{\sigma_i})\cap\gimel(
P_{\sigma_{i+2}})|$, we have $\mfv_{i+1}= \mfv_{j+1}$. Thus $i=j$.
Then all $\sigma_i$ are distinct,  $C$ is an $n$-cycle  and $C$ is
$\langle x,y\rangle$-symmetric. It implies $X_{[C]}^{[C]}\cong
D_{2n}$. Hence $X_{[C^g]}^{[C^g]}\cong D_{2n}$ for any $g\in X$.

Set $\mcE=\{C^x\div x\in X\}$. Then $\mcE$ is an $X$-orbit of
$n$-cycles of $\Sigma$. Since $C$ is $X_{[C]}$-symmetric, $C$ is
$(X_{[C]},3)$-arc-transitive. Recall that the $3$-arc
$(\sigma_{i-1}, \sigma_i, \sigma_{i+1},  \sigma_{i+2}) $ of $C$ is
contained in $\Delta$. It follows that $\Delta=\cup_{C\in
\mcE}Arc_3(C)$.

It is easily shown that $X_{[\mfC]}$ is a
 subgroup of $X_{[C]}$, and so $|\mcE|=|X:X_{[C]}|\le |X:X_{[\mfC]}|=m$.
Suppose that $X_{[\mfC]}$ is a proper subgroup of $X_{[C]}$. Then
there is some $z\in X_{[C]}$ with $C^z=C$ but $\mfC^z\ne \mfC$.
Noting that $\mfC$ and $\mfC^z$ are distinct connected component of
$\gimel$, we have $V(\mfC)\cap V(\mfC^z)=\emptyset$. Since $C^z=C$,
there exist $i$, $j$ and $l$ with $\sigma_1=\sigma_i^z$,
$\sigma_2=\sigma_j^z$ and $\sigma_3=\sigma_l^z$. Then
$\mfv_i^z=[\tau_i^z,\sigma_1,\delta_i^z]\in P_{\sigma_1}$,
$\mfv_j^z=[\tau_j^z,\sigma_2,\delta_j^z]\in P_{\sigma_2}$ and
$\mfv_l^z=[\tau_l^z,\sigma_3,\delta_l^z]\in P_{\sigma_3}$. Since
$(\sigma_1, \sigma_2,\sigma_3)$ is a $2$-arc of $C$, we know
$(\sigma_i,\sigma_j,\sigma_l)$ is also a $2$-arc of $C$. It follows
that $i-j\equiv j-l\equiv \pm 1\, (\mod n)$. Then
$[\mfv_i,\mfv_j,\mfv_l]$ is a $2$-path of $\mfC$, and so
$[\mfv_i^z,\mfv_j^z,\mfv_l^z]$ is a $2$-path of $\mfC^z$. Thus
$\mfv_2,\,\mfv_j^z\in P_{\sigma_2}\cap \gimel(
P_{\sigma_1})\cap\gimel( P_{\sigma_3})$. Since $V(\mfC)\cap
V(\mfC^z)=\emptyset$, we have $\mfv_2\ne \mfv_j^z$, which
contradicts  $\lambda=1$. $X_{[\mfC]}=X_{[C]}$ and so
$|\mcE|=|X:X_{[C]}|=|X:X_{[\mfC]}|=m$.

Since $\Sigma$ is $(X,2)$-arc-transitive, every $2$-path is
contained in some $n$-cycle in $\mcE$. Then
$mn=|Path_2(\Sigma)=|\cup_{C\in\mcE}Path_2(C)|\le
\sum_{C\in\mcE}|Path_2(C)|=mn$. It follows that every  $2$-path of
$\Sigma$ is contained in a unique member of $\mathcal {E}$. Thus
either $girth(\Sigma)=3$ and $\Sigma\cong K_{\mbv+1}$, or $n\ge
girth(\Sigma)\ge 4$ and $\Sigma$ is a near $n$-gonal graph with
respect to $\mathcal {E}$. \qed

 The following result follows
from Theorem~\ref{thm_evenvalency_selfpaired} and
\ref{thm_nearngonal-graph}.

\begin{cor}\label{cor_evenregular_ngonal}
Any connected $(X,2)$-arc-regular graph with even valency and girth
no less than $4$ is a near $n$-gonal graph for some integer $n \geq
4$.
\end{cor}

\section{On tetravalent symmetric graphs}
Let $\Sigma$ be a regular graph with valency four. Recall that
$J(\Sigma)$ is the set of  pairs $([\tau', \tau, \tau''], [\sigma',
\sigma, \sigma''])$ of $2$-paths of $\Sigma$ such that $\sigma \in
\Sigma(\tau) \setminus \{\tau', \tau''\}$, $\tau \in \Sigma(\sigma)
\setminus \{\sigma', \sigma''\}$. For an arbitrary $3$-arc
$\bbalpha:=(\tau_1, \tau, \sigma, \sigma_1)$ of $(\Sigma)$, let
$J_\bbalpha$ be the pair $([\tau_2, \tau, \tau_3]$, $[\sigma_2,
\sigma, \sigma_3])$ of $2$-paths of $\Sigma$, where  $\Sigma(\tau) =
\{\sigma, \tau_1, \tau_2, \tau_3\}$ and $\Sigma(\sigma) = \{\tau,
\sigma_1, \sigma_2, \sigma_3\}$. Then $J_\bbalpha \in J(\Sigma)$.
For any subset $\Delta$ of $Arc_3(\Sigma)$, we set $J({\Delta}) :=
\{J_\bbalpha\div \bbalpha \in \Delta\}$. It is easily shown that
$\Delta$ is a self-paired $X$-orbit on $Arc_3(\Sigma)$ if and only
if $J(\Delta)$ is a self-paired $X$-orbit on $J(\Sigma)$.

\begin{thm}\label{thm_tetravalent_self-paired}
Let $\Sigma$ be a connected $(X,2)$-arc-transitive graph of valency
$4$. If  $\Delta$ is a self-paired $X$-orbit on $Arc_3(\Sigma)$,
then  $\gimel(\Sigma, \Delta) \cong \Psi(\Sigma, J(\Delta))$.
\end{thm}
\proof Define $\psi: Path_2(\Sigma)\rightarrow
Path_2(\Sigma);\,[\tau_1, \tau, \tau_2]\mapsto [\tau_3, \tau,
\tau_4]$, where $\{\tau_3, \tau_4\} = \Sigma(\tau)\setminus
\{\tau_1, \tau_2\}$. It is easy to check that $\psi$ is an
isomorphism from  $\gimel(\Sigma, \Delta)$ to $\Psi(\Sigma,
J(\Delta))$. \qed

The main aim of this section is to give a characterization of
tetravalent $(X,2)$-arc-transitive graphs. The following simple
lemma is useful.
\begin{lem}\label{connectedness} Let $\Gamma$ be an $X$-symmetric
graph with an $X$-invariant partition $\mcB$ such that $\Gamma_\mcB$
is connected and $(X,2)$-arc-transitive. Let $B\in \mcB$ and
$C,\,D\in \Gamma_\mcB(B)$ with $C\ne D$. If $\Gamma[B,C]$ is
connected and $\Gamma(C)\cap B\cap \Gamma(D)\ne \emptyset$, then
$\Gamma$ must be connected.
\end{lem}
\proof It suffices to show that there is a path in $\Gamma$ between
any two different vertices $\mfv$ and $\mfu$ of $\Gamma$. Since
$\Gamma_\mcB$ is $(X,2)$-arc-transitive, $\Gamma[B,C]$ is
independent of the choices of $B$ and $C\in \Gamma_\mcB(B)$ up to
isomorphism; and $|\Gamma(C)\cap B\cap \Gamma(D)|$ is independent of
the choices of $B$ and $C,\, D\in \Gamma_\mcB(B)$ (with $C\ne D$).

Assume first $\mfv,\, \mfu\in B$. Without loss of generality, we
assume $\mfv\in \Gamma(C)\cap B\cap \Gamma(D)$. If $\mfu\in
\Gamma(C)\cap B$, then there a path in $\Gamma$ between $\mfv$ and
$\mfu$ as $\Gamma[B,C]$ is connected. So we assume $\mfu\not\in
\Gamma(C)\cap B$. Take $E\in \Gamma_\mcB(\mfu)$. Then $E\in
\Gamma_\mcB(B)$, $\mfu\in B\cap \Gamma(E)$ and $|\Gamma(C)\cap B\cap
\Gamma(E)|=|\Gamma(C)\cap B\cap \Gamma(D)|>0$. Let $\mfw\in
\Gamma(C)\cap B\cap \Gamma(E)$. Then either $\mfv=\mfw$ or there is
a path between $\mfv$ and $\mfw$, and there is a path between $\mfw$
and $\mfu$. Thus there is a path between $\mfv$ and $\mfu$.

Now let $\mfv\in B$ and $\mfu\in B'$ with $B\ne B'$. Since
$\Gamma_\mcB$ is connected, there is a path $[B=B_1,\ldots,
B_l=B']$. Let $\mfu_l'\in B_l$ and $\mfu_{l-1}\in B_{l-1}$ such that
$\{\mfu_{l-1}, \mfu_l'\}\in E(\Gamma)$. Thus there is a path between
$\mfu_{l-1}$ and $\mfu$. Then induction on $l$ implies that there is
a path between $\mfv$ and $\mfu$. \qed

Now we are ready to state and prove the main result of this section.
\begin{thm}\label{thm_tetravalent}
Let $\Sigma$ be a connected $(X,2)$-arc-transitive graph with
valency $\mathbbm{v} = 4$, where $X$ acts faithfully on $V(\Sigma)$.
Then $\Sigma$ has a self-paired $X$-orbit $\Delta$ on
$Arc_3(\Sigma)$. Set
 $\gimel := \gimel(\Sigma, \Delta)$, $\Xi := \Xi(\Sigma, \Delta)$, $e :=
|E(\Sigma)|$, $\mu := |V(\Sigma)|$. Let $(\tau, \sigma) \in
Arc(\Sigma)$. Then one of the following cases occurs.
\begin{enumerate}
\item[(a)]
$\gimel[P_{\tau}, P_{\sigma}] \cong \Xi[A_{\tau}, A_{\sigma}] \cong
3K_2$, $\gimel(\Sigma, \Delta) \cong mC_n$, $val(\Xi) = 3$, and
there exists an $X$-orbit $\mathcal {E}$ of $n$-cycles of $\Sigma$
with $|\mathcal {E}| = m$, such that $\Delta=\cup_{C\in
\mcE}Arc_3(C)$, $X_{[C]}^{[C]} \cong D_{2n}$ for each $C\in \mathcal
{E}$, where $m \geq 6$ and $n \geq girth(\Sigma)$ with $mn = 3e =
6\mu$. Moreover, either $\Sigma \cong K_{5}$ or $\Sigma$ is a near
$n$-gonal graph with respect to $\mathcal {E}$; and, either
\begin{enumerate}
\item[(a.1)]
$X_{P_{\tau}} =X_{A_{\tau}}=X_\tau \cong A_4$, both $\gimel$ and
$\Xi$ are $(X,1)$-arc-regular and $\Sigma$ is $(X,2)$-arc-regular;
or
\item[(a.2)]
$X_\tau =X_{P_{\tau}} = X_{A_{\tau}}  \cong S_4$, $\gimel$ is  not
$(X,1)$-arc-regular, $\Xi$ is $(X,2)$-arc-regular.
\end{enumerate}
\item[(b)]
$\gimel[P_{\tau}, P_{\sigma}] \cong \Xi[A_{\tau}, A_{\sigma}] \cong
K_{3,3} - 3K_2$, $val(\gimel) = 4$, $val(\Xi) = 6$, $X_{P_{\tau}}=
X_{A_{\tau}}=X_\tau\cong S_4$, both $\gimel$ and $\Xi$ are connected
and $(X,1)$-arc-regular. Further, $\Delta' := Arc_3(\Sigma)
\setminus \Delta$ is a self-paired $X$-orbit on $Arc_3(\Sigma)$, and
there exists an $X$-orbit $\mathcal {E}$ of $n$-cycles of $\Sigma$
with $|\mathcal {E}| = m$, such that $\Delta' = \cup_{C \in \mathcal
{E}}Arc_3(C)$, $X_{[C]}^{[C]} \cong D_{2n}$ for each $C \in \mathcal
{E}$, where $m \geq 6$ and $n \geq girth(\Sigma) \geq 3$ with $mn =
3e = 6\mu$. Moreover, either $\Sigma \cong K_5$ or $\Sigma$ is a
near $n$-gonal graph with respect to $\mathcal {E}$.
\item[(c)]
$\gimel[P_{\tau}, P_{\sigma}] \cong \Xi[A_{\tau}, A_{\sigma}] \cong
K_{3,3}$, $val(\gimel) = 6$, $val(\Xi) = 9$, both $\gimel$ and $\Xi$
are connected and $(X,1)$-transitive, and $\Sigma$ is
$(X,3)$-arc-transitive.
\end{enumerate}
\end{thm}
\proof By Theorem \ref{thm_evenvalency_selfpaired}, $\Sigma$ has a
self-paired $X$-orbit $\Delta$ on $Arc_3(\Sigma)$. Then, by
Proposition~\ref{prop_U-graph}, $\gimel := \gimel(\Sigma, \Delta)$
is $X$-symmetric and admits an $X$-invariant partition
$\mcP:=\{P_\sigma\div \sigma\in V(\Sigma)\}$ with $\Sigma\cong
\gimel_\mcP$, and by Proposition~\ref{prop_3-arc_graph}, $\Xi :=
\Xi(\Sigma, \Delta)$ is $X$-symmetric and  admits an $X$-invariant
partition $\mathcal{A}:=\{A_\sigma\div \sigma\in V(\Sigma)\}$ with
$\Sigma\cong \Xi_\mathcal{A}$. Let $\ell_i := \ell_i(\Delta)$, $i =
1, 2, \ldots, t$, be defined as in Section 5. Then $t\le 3$ as
$val(\Sigma)=4$.

Let $(\tau,\sigma)\in Arc(\Sigma)$. Then there is a $3$-arc
$(\tau_1,\tau,\sigma,\sigma_1)\in \Delta$ as $\Sigma$ is
$X$-symmetric. It follows that
$\{[\tau_1,\tau,\sigma],[\tau,\sigma,\sigma_1]\}$ is an edge of
$\gimel[P_{\tau}, P_{\sigma}]$, and that
$\{(\tau,\tau_1),(\sigma,\sigma_1)\}$ is an edge of $\Xi[A_{\tau},
A_{\sigma}]$. It is easily shown that $X_{(\tau,\sigma)}=X_\tau\cap
X_\sigma=X_{P_\tau}\cap X_{P_\sigma}$ acts transitively on the edges
of $\gimel[P_{\tau}, P_{\sigma}]$. It implies that the stabilizer
$(X_{(\tau,\sigma)})_{[\tau_1,\tau,\sigma]}=X_{(\tau_1,\tau,\sigma)}$
acts transitively on the neighborhood of $[\tau_1,\tau,\sigma]$ in
$\gimel[P_{\tau}, P_{\sigma}]$. Then the valency of
$\gimel[P_{\tau}, P_{\sigma}]$ equals to
$|X_{(\tau_1,\tau,\sigma)}:(X_{(\tau_1,\tau,\sigma)})_{[\tau,\sigma,\sigma_1]}|
=|X_{(\tau_1,\tau,\sigma)}:X_{(\tau_1,\tau,\sigma,\sigma_1)}|=\ell_1$.
Further, since $\Sigma$ is $(X,2)$-arc-transitive,
$X_{(\tau,\sigma)}$ is transitive on
$\Sigma(\tau)\setminus\{\sigma\}:=\{\tau_1,\tau_2,\tau_3\}$ and on
$\Sigma(\sigma)\setminus\{\tau\}:=\{\sigma_1,\sigma_2,\sigma_3\}$.
Thus $V(\gimel[P_{\tau}, P_{\sigma}])=\{[\tau_i,\tau,\sigma]\div
i=1,2,3\}\cup \{[\tau,\sigma,\sigma_i]\div i=1,2,3\}$. A similar
argument leads to $V(\Xi[A_{\tau}, A_{\sigma}])=\{(\tau,\tau_i)\div
i=1,2,3\}\cup \{(\sigma,\sigma_i)\div i=1,2,3\}$. It is easy to
check that $[\tau_i,\tau,\sigma]\mapsto
(\tau,\tau_i),\,[\tau,\sigma,\sigma_i]\mapsto (\sigma,\sigma_i)$
gives an isomorphism from $\gimel[P_{\tau}, P_{\sigma}]$ to
$\Xi[A_{\tau}, A_{\sigma}]$. Further, $\gimel[P_{\tau},
P_{\sigma}]\cong 3K_2$, $K_3-3K_2$ or $K_{3,3}$ according to
$\ell_1=1$, $2$ or $3$, respectively. By \cite[Theorem
4.3]{FSGTTQ(II)}, $2=r_\mcP :=|\gimel_\mcP([\tau_1,\tau,\sigma])|$
for any $[\tau_1,\tau,\sigma]\in V(\gimel)$. Then
$val(\gimel)=r_\mcP \ell_1=2\ell_1$. By Lemma~\ref{3-arc-graph},
$val(\Xi)=r_\mathcal{A} \ell_1=3\ell_1$. Since $\Sigma$ is
$(X,2)$-arc-transitive, $X_\tau^{\Sigma(\tau)}\cong A_4$ or $S_4$.
It is easy to see $X_\tau=X_{P_\tau}=X_{A_\tau}$,
$(X_\tau)_{(\Sigma(\tau))}=X_{(P_\tau)}=X_{(A_\tau)}$ and hence
$X_\tau^{\Sigma(\tau)}\cong
X_{P_\tau}^{P_\tau}=X_{A_\tau}^{A_\tau}$. We treat the following
three separate cases.

{\bf Case 1}. $\ell_1=1$. Then $val(\gimel)=2$, $val(\Xi)=3$ and
$\bfl(\Delta)=(1,1,1)$ or $(1,2)$ in this case.  By Theorem
\ref{thm_nearngonal-graph}, the part of (a) prior to (a.1) holds.
Again by Theorem~\ref{thm_nearngonal-graph}, $X_\tau$ acts
faithfully on $\Sigma(\tau)$, and hence $X_\tau^{\Sigma(\tau)}\cong
X_\tau$.

Assume first $\bfl(\Delta)=(1,1,1)$. Then $X_{(\tau_1, \tau,
\sigma)}\le X_{(\tau, \sigma,\sigma_1)}$ for any $2$-arc $(\tau_1,
\tau, \sigma)$ of $\Sigma$ and $\tau\ne\sigma_1\in \Sigma(\sigma)$.
Since $\Sigma$ is $(X,2)$-arc-transitive, the stabilizers of any two
$2$-arcs of $\Sigma$ are conjugate in $X$, in particular, they has
the same order. Thus $X_{(\tau_1, \tau, \sigma)}=X_{(\tau,
\sigma,\sigma_1)}$. Since $\Sigma$ is connected, $X_{(\tau_1, \tau,
\sigma)}=X_{(\tau_1', \tau', \sigma')}$ for an arbitrary $2$-arc
$(\tau_1', \tau', \sigma')$ of $\Sigma$. Hence $X_{(\tau_1, \tau,
\sigma)}=1$ as $X$ is faithful on $V(\Sigma)$. Then $\Sigma$ is
$(X,2)$-arc-regular. It implies $X_{P_\tau}=X_{A_\tau}=X_\tau\cong
A_4$. Then (a.1) follows from calculating the numbers of arcs or
$2$-arcs of $\gimel$, $\Xi$ and $\Sigma$.

Now let $\bfl(\Delta)=(1,2)$. Then  $X_{(\tau_1, \tau, \sigma)}$
acts transitively on $\Sigma(\sigma)\setminus \{\tau,\sigma_1\}$ for
$(\tau_1, \tau, \sigma,\sigma_1)\in \Delta $. Thus $X_{(\tau_1,
\tau, \sigma)}\ne 1$, and $\Sigma$ is not $(X,2)$-arc-regular.
Recall that $X_\tau$ acts faithfully on $\Sigma(\tau)$. It implies
$X_{P_{\tau}}  =X_{A_{\tau}}  =X_\tau \cong S_4$. Since $\ell_1=1$,
we have $X_{(\tau_1, \tau, \sigma,\sigma_1)}=X_{(\tau_1, \tau,
\sigma)}\ne 1$. It implies that $\gimel$ is not $(X, 1)$-arc
regular.

Let $(\tau,\sigma)\in V(\Xi)=Arc(\Sigma)$. Set
$\Xi((\tau,\sigma))=\{(\sigma_1,\delta_1),(\sigma_2,\delta_2),(\sigma_3,\delta_3)\}$,
the neighborhood of $(\tau,\sigma)$ in $\Xi$. Then
$\Sigma(\tau)=\{\sigma,\sigma_1,\sigma_2,\sigma_3\}$ and $(\sigma,
\tau, \sigma_i,\delta_i)\in \Delta$, $i=1,\,2,\,3$. It follows from
$\ell_1=1$ that $\sigma_i^x=\sigma_j$ implies $\delta_i^x=\delta_j$
for $x\in X_{(\tau,\sigma)}$ and $1\le i,\,j\le 3$. Then
$X_{(\tau,\sigma)}$ fixes $\Xi((\tau,\sigma))$ setwise. Since
$X_\tau \cong S_4$, we conclude that the permutation group induced
by $X_{(\tau,\sigma)}$ on $\Sigma(\tau)\setminus \{\sigma\}$ is
isomorphic to $S_3$, which is $2$-transitive on
$\Sigma(\tau)\setminus \{\sigma\}$. Thus  $X_{(\tau,\sigma)}$ acts
$2$-transitively on $\Xi((\tau,\sigma))$. It follows that $\Xi$ is
$(X,2)$-arc-transitive. Further, checking the number of the $2$-arcs
of $\Xi$ implies that $\Xi$ is $(X,2)$-arc-regular. This complete
the proof of (a).

{\bf Case 2}. $\ell_1=2$. In this case, $val(\gimel) = 2\ell_1 = 4$,
$val(\Xi) = 3\ell_1 = 6$ and $\gimel[P_{\tau}, P_{\sigma}] \cong
\Xi[A_{\tau}, A_{\sigma}] \cong K_{3,3} - 3K_2$. By
Lemma~\ref{connectedness}, both $\gimel$ and $\Xi$ are connected.

Now we shall show that $X_\tau$ acts faithfully on the neighborhood
$\Sigma(\tau)$ of $\tau$ in $\Sigma$, by a similar argument as in
the first paragraph of the proof of
Theorem~\ref{thm_nearngonal-graph}. Since $\Sigma$ is $(X,2)$-arc
transitive, every $2$-arc of $\Sigma$ lies in a member of $\Delta$.
Let $(\tau,\sigma)$ be an arbitrary arc of $\Sigma$. Since
$\ell_1(\Delta)=2$ and $\Delta$ is a self-paired $X$-orbit, we
conclude that, for any $\tau_1\in \Sigma(\tau)\setminus\{\sigma\}$
there is a unique $\sigma_1\in \Sigma(\sigma)\setminus\{\tau\}$ such
that $(\tau_1,\tau,\sigma,\sigma_1)\not\in \Delta$,
$X_{(\tau_1,\tau,\sigma)}=X_{(\tau,\sigma,\sigma_1)}$, and that
$(\tau_1',\tau,\sigma,\sigma_1)\not\in \Delta$ implies
$\tau'_1=\tau_1$. Then $(X_\tau)_{(\Sigma(\tau))}=\cap_{\tau_1\in
\Sigma(\tau)\setminus\{\sigma\}}X_{(\tau_1,\tau,\sigma)}=\cap_{\sigma_1\in
\Sigma(\sigma)\setminus\{\tau\}}X_{(\tau,\sigma,\sigma_1)}=(X_\sigma)_{(\Sigma(\sigma))}$.
It follows from the connectedness of $\Sigma$ that
$(X_\tau)_{(\Sigma(\tau))}$ fixes every vertex of $\Sigma$. Thus
$(X_\tau)_{(\Sigma(\tau))}=1$ and $X_\tau$ is faithful on
$\Sigma(\tau)$.

For a $2$-arc $(\tau_1,\tau,\sigma)$ of $\Sigma$, since
$\ell_1(\Delta)=2$, there is $\sigma_2,\sigma_3\in \Sigma(\sigma)$
such that $(\tau_1,\tau,\sigma,\sigma_2)\in \Delta$ and
$(\tau_1,\tau,\sigma,\sigma_3) \in \Delta$. Since $\Delta$ is an
$X$-orbit, $X_{(\tau_1,\tau,\sigma)}$ acts transitively on
$\{\sigma_2,\sigma_3\}$. In particular, $X_{(\tau_1,\tau,\sigma)}\ne
1$. Thus we have $X_{P_{\tau}}  =X_{A_{\tau}}  =X_\tau\cong S_4$.
Further, $|X|=|V(\Sigma)||X_\tau|=24\mu=|Arc(\gimel)|=|Arc(\Xi)|$,
so both $\gimel$ and $\Xi$ are $(X,1)$-arc-regular.

Set $\Delta' = Arc_3(\Sigma) \setminus \Delta$. Then $\Delta'$ is
self-paired and $X$-invariant.  For any two $3$-arcs $(\tau_1, \tau,
\sigma, \sigma_1)$ and $(\tau_1', \tau', \sigma', \sigma_1')$ of
$\Sigma$ in $\Delta'$, since $\Sigma$ is $(X,2)$-arc-transitive,
there exists some $x \in X$ such that $(\tau_1', \tau', \sigma')^x =
(\tau_1, \tau, \sigma)$. Then $(\tau_1, \tau, \sigma, \sigma_1'^x) =
(\tau_1', \tau', \sigma', \sigma_1')^x \in \Delta'$. By the argument
in the second paragraph of this case,  $\sigma_1'^x = \sigma_1$,
that is, $(\tau_1', \tau', \sigma', \sigma_1')^x = (\tau_1, \tau,
\sigma, \sigma_1)$. It follows that $\Delta'$ is $X$-transitive and
$\ell_1(\Delta')=1$. Thus (b) holds by
Theorem~\ref{thm_nearngonal-graph}.

{\bf Case 3}.  $\ell_1(\Delta) = 3$. Then $val(\gimel) = 2\ell_1 =
6$, $val(\Xi) = 3\ell_1 = 9$ and $\gimel[P_{\tau}, P_{\sigma}] \cong
\Xi[A_{\tau}, A_{\sigma}] \cong K_{3,3}$. It follows from
\cite[Theorem 2]{ACOFSGWTTQ} that  $\Sigma$ is $(X,3)$-arc
transitive. By Lemma~\ref{connectedness}, both $\gimel$ and $\Xi$
are connected. Note that $\Sigma\cong \gimel_\mcP$ is of valency
four. Let $\sigma,\, \tau$ and $\delta$ be three distinct vertices
of $\Sigma$ such that $P_\sigma,\,P_\delta\in \gimel_\mcP(P_\tau)$.
Then there exist $\mfv\in P_\tau$, $\mfu_1,\mfu_2\in P_\sigma$ and
$\mathfrak{w}\in P_\delta$ such that $(\mfu_1,\mfv, \mfu_2)$ and
$(\mathfrak{w},\mfv, \mfu_1)$ are $2$-arcs of $\gimel$. Since $\mcP$
is $X$-invariant, there is no $x\in X$ with $(\mfu_1,\mfv,
\mfu_2)^x=(\mathfrak{w},\mfv, \mfu_1)$. Thus $\gimel$ is not
$(X,2)$-arc-transitive, and so it is $(X,1)$-transitive. A similar
argument implies that $\Xi$ is $(X,1)$-transitive. Hence (c) holds.
\qed

\begin{cor}\label{4-val-near}
Let $\Sigma$ be a connected tetravalent $(X,2)$-transitive graph.
Then either $\Sigma\cong K_5$, or $\Sigma$ is a near $n$-gonal graph
for some integer $n \geq 4$.
\end{cor}
At the end of this section we give several examples, which indicate
that there exist certain graphs satisfying each case listed in
Theorem~\ref{thm_tetravalent}.

\begin{ex}\label{ex_tetravalent_A1}
{\rm Let $X = PSL(2,p)$, where $p$ is a prime such that $5\ne p
\equiv \pm 3 \,({\rm mod}\, 8)$. Then by \cite{FSTGOOO}, there exist
$H < X$ and an involution $z\in X$  such that $H \cong A_4$, $P = H
\cap H^z =: \langle h \rangle \cong Z_3$, $z \in N_X(P)$ and $h^z =
h^{-1}$. Moreover, $\Sigma := Cos(X, H, HzH) \ncong K_5$ is a
tetravalent $(X,2)$-arc-regular graph and $Aut(\Sigma) = X$. Set $H
= P\cup Pg\cup Pg_2\cup Pg_3$. Let $\Delta=\{(Hzgx, Hx, Hzx,
Hzgzx)\div x\in X\}$. Then $\Delta$ is a self-paired $X$-orbit on
$Arc_3(\Sigma)$ with $\bfl(\Delta) = (1, 1, 1)$. }
\end{ex}

\begin{ex}\label{ex_tetravalent_A2}
{\rm Let  $X = PSL(2,p)$ for a prime $p\ge 11$ with $p \equiv \pm
1\,({\rm mod}\, 8)$. Let $S_4\cong H < X$. Then by \cite[Lemma
4.1]{OPPGWSSATOG}, there exists an involution $z \in X \setminus H$
such that $N_{X}(P) = P \times \langle z\rangle$, where $P = H \cap
H^z \cong S_3$. Further, $\Sigma = Cos(X, H, HzH)$ is a tetravalent
$(X,2)$-transitive graph with $Aut(\Sigma) = X$. Set $H = P\cup
Pg\cup Pg_2\cup Pg_3$. Then $\Delta=\{(Hzgx, Hx, Hzx, Hzgzx)\div
x\in X\}$ is a self-paired $X$-orbit on $Arc_3(\Sigma)$ with
$\bfl(\Delta) = (1,2)$. }
\end{ex}

\begin{ex}\label{ex_tetravalent_A2B}
{\rm Let $\Sigma =K_{5,5}-5K_2$ with vertex set $\{i,i'\div 1\le
i\le 5\}$. For  $g\in S_5$, define ${\bar g}: i\mapsto g(i),\,
i'\mapsto g(i)'$. Let $z: i\leftrightarrow i'$. Set $X=\langle {\bar
g}, z\div g\in S_5\rangle$. Then $\Sigma$ is $(X,2)$-transitive.
Then both $\Delta_1:=\{(1,2',3,1')^x\div x\in X\}$ and
$\Delta_2:=\{(1,2',3,4')^x\div x\in X\}$ are self-paired with
$\bfl(\Delta_1) = (1, 2)$ and $\bfl(\Delta_2) = (2,1)$. }
\end{ex}

%
%
\section{Heptavalent graphs with $X_\tau^{\Sigma(\tau)} \cong PSL(3,2)$}
\begin{thm}\label{thm_heptavalent}
Let $\Sigma$ be an $(X,2)$-arc-transitive graph of valency $7$ with
$X_\tau^{\Sigma(\tau)} \cong PSL(3,2)$ for $\tau\in V(\Sigma)$. Then
there exists a self-paired $X$-symmetric orbit $\Theta$ on $D\mcS
t^{3}(\Sigma)$. Let $\Pi = \Pi(\Sigma, \Theta)$ and $\mcS = \mcS
t(\Theta)$. Then, for $\sigma \in \Sigma(\tau)$, one of the
following cases occurs.
\begin{enumerate}
\item[(1)] $\Pi[\mcS_{\tau}, \mcS_{\sigma}] \cong 3K_2$, and
$\Pi$ is a trivalent $(X,2)$-arc-transitive graph;
\item[(2)] $\Pi[\mcS_{\tau}, \mcS_{\sigma}] \cong K_{3,3}-3K_2$,
$val(\Pi)=6$, $\Pi$ is connected and $(X,1)$-transitive;
\item[(3)] $\Pi[\mcS_{\tau}, \mcS_{\sigma}] \cong K_{3,3}$, $val(\Pi)=9$  and
$\Pi$ is  connected and $(X,1)$-transitive.
\end{enumerate}
\end{thm}
\proof Let $\tau\in V(\Sigma)$. Since $X_\tau^{\Sigma(\tau)} \cong
PSL(3,2)$, we may identify $\Sigma(\tau)$ with the point set of
seven-point plane $PG(2,2)$, which is an $X_{\tau}$-flag-transitive
$1$-$(7,3,3)$ design with multiplicity $1$. By Theorem
\ref{thm_design_odd}, there exists a self-paired $X$-symmetric orbit
$\Theta$ on $D\mcS t^3(\Sigma)$. Set $\mcS = \mcS t(\Theta)$ and
$\Pi = \Pi(\Sigma, \Theta)$. Then, by
Theorem~\ref{thm_doublestar_graph}, $\Pi$ is $X$-symmetric and
$\Pi_\mcB\cong \Sigma$, where $\mcB=\{\mcS_\tau\div \tau\in
V(\Sigma)\}$ and $\mcS_\tau=\{\mfs\in \mcS\div \mfs=\mfs(\tau,S),
S\subseteq\Sigma(\tau), |S|=3\}$. Further, for $\mcS_\tau\in \mcB$,
we have $X_\tau=X_{\mcS_\tau}$ and  $\mcD(\mcS_\tau)\cong
\mbD^*(\tau)\cong PG(2,2)$. (See Section 3 for the definition of
$\mbD(\tau)$.) In particular, for $\sigma\in \Sigma(\tau)$,
$|\mcS_\tau\cap\Pi(\mcS_\sigma)|=3$; thus the bipartite graph
$\Pi[\mcS_{\tau}, \mcS_{\sigma}]$ is isomorphic to one of $3K_2$,
$K_{3,3}-3K_2$ and $K_{3,3}$ as $X_\tau\cap X_\sigma$ acts
transitively on the edges of $\Pi[\mcS_{\tau}, \mcS_{\sigma}]$.
Moreover, noting that, any pair of distinct lines of $PG(2,2)$
intersect a unique point and any pair of distinct points determine a
unique line, it follows that $\lambda:=|\Pi(\mcS_{\sigma})\cap
\mcS_{\tau}\cap \Pi(\mcS_{\delta})|=1$ for $\sigma,\delta\in
\Sigma(\tau)$ with $\sigma\ne \delta$. Then by
Lemma~\ref{connectedness}, $\Pi$ is connected if $\Pi[\mcS_{\tau},
\mcS_{\sigma}]\cong K_{3,3}-3K_2$ or $K_{3,3}$. Note that each point
of $\mcD(\mcS_\tau)$ belongs to three blocks. It follows that $\Pi$
is of valency $3\ell$, where $\ell$ is the valency of
$\Pi[\mcS_{\tau}, \mcS_{\sigma}]$.

Assume first that $\Pi[\mcS_{\tau}, \mcS_{\sigma}]\cong 3K_2$. Then
$val(\Pi)=3$. Let $\mfs\in \mcS_\tau$, and
$\Pi(\mfs)=\{\mfs_1,\mfs_2,\mfs_3\}$ with $\mfs_i\in \mcS_{\tau_i}$
for $i=1,2,3$. Then $\tau_1$, $\tau_2$ and $\tau_3$ are distinct
vertices of $\Sigma$. Recall $\mcD(\mcS_\tau)\cong \mbD^*(\tau)\cong
PG(2,2)$. Then we may identify $\mfs$ with a line $L$ of $PG(2,2)$,
and $\mcS_{\tau_i}$ with the points in this line. Then
$(X_\tau^{\Sigma(\tau)})_\mfs\cong S_4$ acts $2$-transitively on
$\{\mcS_{\tau_i}\div i=1,2,3\}$. It implies that
$(X_\tau)_\mfs=X_\mfs$ acts $2$-transitively (and unfaithfully) on
$\{\mfs_1,\mfs_2,\mfs_3\}$. Thus $\Pi$ is $(X,2)$-arc-transitive,
and (1) holds.

Now let $\Pi[\mcS_{\tau}, \mcS_{\sigma}]\cong K_{3,3}-3K_2$ or
$K_{3,3}$. Then $\Pi$ has two $2$-arcs, say $(\mfv,\mfu,\mfw)$ and
$(\mfv',\mfu',\mfw')$, such that $\mfv,\mfv',\mfw\in \mcS_\tau$,
$\mfu,\mfu'\in \mcS_\sigma$ and $\mfw'\in \mcS_\delta$ for distinct
$\tau$, $\sigma$ and $\delta$. Noting $\mcB$ is $X$-invariant, there
is no $x\in X$ maps $(\mfv,\mfu,\mfw)$ to $(\mfv',\mfu',\mfw')$.
Thus $\Pi$ is not $(X,2)$-arc-transitive. Then (2) and (3) hold.
\qed

The following  examples indicate that there exist certain graphs
satisfying each case listed in Theorem~\ref{thm_heptavalent}.

\begin{ex}\label{ex_heptavalent_AB}
{\rm Let $\Sigma$ be the complete graph on vectors of
$\mathbb{F}^3$, where $\mathbb{F} = \{0, 1\}$ is a binary field.
Then the $3$-dimensional affine group $X:=AGL(3,2)$ is a subgroup of
the automorphism group $Aut(\Sigma)\cong S_8$ of $\Sigma$. Set
$\bfv_0=(0,0,0)$, $\bfv_1=(1,0,0)$, $\bfv_2=(0,1,0)$,
$\bfv_3=(1,1,0)$, $\bfv_4=(0,0,1)$, $\bfv_5=(1,0,1)$,
$\bfv_6=(0,1,1)$ and $\bfv_7=(1,1,1)$. Then $X_{\bfv_0}=GL(3,2)\cong
PSL(3,2)$ is $2$-transitive on $\{\bfv_i\div i=1,2,\ldots, 7 \}$.
Hence $\Sigma$ is $(X,2)$-arc-transitive. We define
$t_1:(a_1,a_2,a_3)\mapsto (a_1,a_2,a_3+1)$ and
$t_2:(a_1,a_2,a_3)\mapsto (a_2,a_1,a_3+1)$  for $(a_1,a_2,a_3)\in
\mathbb{F}^3$, respectively. Then $t_1,\,t_2\in X$ with
$t_1^2=t_2^2=1$. Let $L= \{\bfv_2, \bfv_4, \bfv_6\}$ and set $\mfl=
\mfs(\bfv_0,L)$. Note that $\{\bfv_0,\bfv_2, \bfv_4, \bfv_6\}$ is  a
subspace of $\mathbb{F}^3$. Then $X_\mfl$  is the stabilizer of this
subspace in $GL(3,2)$. Thus
$$
\begin{array}{l}
X_\mfl=\left\{ {\small \left[\begin{array}{ccc}
1 & e & f \\
0 & a & b \\
0 & c & d
\end{array} \right]}\left|\begin{array}{l} a,b,c,d,e,f\in\mathbb{F}\\
ad-bc=1\end{array}\right.\right\},\\ (X_\mfl)_{\bfv_4}=\left\{
{\small \left[\begin{array}{ccc}
1 & e & f \\
0 & 1 & b \\
0 & 0 & 1
\end{array} \right]}\left|\begin{array}{l} b,e,f\in\mathbb{F} \end{array}\right.\right\}.
\end{array}
$$

Let $\mfr_i= \mfs(\bfv_4,L^{t_i})$ for $i=1$ and $2$. Then
$L^{t_1}=\{\bfv_0,\bfv_2, \bfv_6\}$, $L^{t_2}=\{\bfv_0,\bfv_1,
\bfv_5\}$, $\mfl^{t_i}=\mfr_i$ and $\mfr_i^{t_i}=\mfl$. Thus
$\Theta_i:=\{(\mfl^x,\mfr_i^x)\div x\in X\}$ is a self-paired
$X$-orbits on $D\mcS t^3(\Sigma)$. Let $\Pi^i=\Pi(\Sigma, \Theta_i)$
and $\Delta_i=\Pi^i[\mcS_{\bfv_0}, \mcS_{\bfv_4}]$ for $i=1$ and
$2$. Note that $X_{\mcS_{\bfv_0}}\cap
X_{\mcS_{\bfv_4}}=X_{\bfv_0}\cap X_{\bfv_4}$ acts transitively on
the edges of $\Delta_i$. It follows that $(X_{\bfv_0}\cap
X_{\bfv_4})_\mfl=(X_\mfl)_{\bfv_4}$ is transitive on the
neighborhood of $\mfl$ in $\Delta_i$. Thus
$val(\Delta_i)=|\{\mfr_i^x\div x\in (X_\mfl)_{\bfv_4}\}|$. If $i=1$,
then $\mfr_1^x=\mfs(\bfv_4,L^{t_1x})=\mfs(\bfv_4,L^{t_1})=\mfr_1$
for $x\in (X_\mfl)_{\bfv_4}$, so $val(\Delta_1)=1$ and Theorem
\ref{thm_heptavalent} (1) occurs. (In fact, $\Pi_1\cong 14K_4$. We
omit the detail.) If $i=2$, then $L^{t_2 x}=L^{t_2}$ or $\{\bfv_0,
\bfv_3,\bfv_7\}$ for $x\in (X_\mfl)_{\bfv_4}$, thus
$val(\Delta_2)=2$ and Theorem \ref{thm_heptavalent} (2) occurs. }
\end{ex}

\begin{ex}\label{ex_heptavalent_C}
{\rm Let $\mathbb{F} = \{0, 1\}$ be a  a binary field. Denote by
${\bf i}$ the non-zero vector of $\mathbb{F}^3$ with coordinate
$(a_1, a_2, a_3)$ such that $i = 4a_1 + 2a_2 + a_3$. Let $\Sigma$ be
the complete bipartite graph with vertex set $\{l{\bf i}\div 1\le
i\le 7\}\cup\{r{\bf i}\div 1\le i\le 7\}$. Then $X:=PSL(3,2) \wr
Z_2$ is a subgroup of $Aut(\Sigma)$, and $\Sigma$ is
$(X,3)$-transitive. Let $L= \{r{\bf 1}, r{\bf 2}, r{\bf 3}\}$ and
$R=\{l{\bf 1}, l{\bf 2}, l{\bf 3}\}$. Set $\mfl=\mfs (l{\bf 1},L)$,
$\mfr=\mfs (r{\bf 1},R)$ and $\Theta_3:= \{(\mfl,\mfr)^x\div x\in
X\}$. Then $\Theta_3$ is a self-paired $X$-symmetric orbit on $D\mcS
t^3(\Sigma)$, and $\Pi := \Pi(\Sigma, \Theta)$ satisfies Theorem
\ref{thm_heptavalent} (3).}
\end{ex}

\section{Proof of Theorem \ref{thm_main}}

Now we are ready to give the proof of of Theorem \ref{thm_main}.

Since $\Gamma$ is $X$-symmetric and $\Gamma_\mcB$ contains at least
one edge, $\Gamma_\mcB$ is $X$-symmetric, that is, $X_B$ is
transitive on $\Gamma_\mcB(B)$ for $B\in \mcB$; further, $B$ is an
independent subset of $V(\Gamma)$.

We first show that each of Theorem \ref{thm_main}{(a)-(d)} implies
the $(X,2)$-arc-transitivity of $\Gamma_\mcB$. It suffices to show
that $X_B$ acts $2$-transitively on $\Gamma_\mcB(B)$ for $B\in
\mcB$. It is trivial for the case (d) as $\Gamma_\mcB(B)$  is the
block set of $\mcD(B)$. In the following we assume one of (a), (b)
and (c) occurs.

Suppose that $m:=m(\mcD(B))\ne 1$. Then $\Gamma_\mcB(B)$ admits an
$X_B$-invariant partition $\mathcal{M}:=\{\mathcal{M}_C\div C\in
\Gamma_\mcB(B)\}$, where $\mathcal{M}_C$ is a set of blocks of
$\mcD(B)$ with the same trace $B\cap \Gamma(C)$ of $C$. Thus
$m=|\mathcal{M}_C|$ is a divisor of $b$. For $\mfv\in B$, it is
easily to see that $C\in \Gamma_\mcB(\mfv)$ yields $D\in
\Gamma_\mcB(\mfv)$ for any $D\in \mathcal{M}_C$. This observation
says that $m=|\mathcal{M}_C|$ is a divisor of
$r:=|\Gamma_\mcB(\mfv)|$. It follows that $(v,b,r)=(6,4,2)$, $m=2=r$
and $|\mathcal{M}|=2$. Set $\mathcal{M}=\{\mathcal{M}_C,
\mathcal{M}_D\}$. Then $\mathcal{T}:=\{B\cap \Gamma(C), B\cap
\Gamma(D)\}$ is an $X_B$-invariant partition of $B$. Let $K$ be the
kernel of $X_B$ acting on $\mathcal{T}$. Then $|X_B:K|=2$ and
$X_{(B)}\le K$. It follows that $X_B^B\cong S_4$ and $K/X_{(B)}\cong
A_4$. Note that $K$ is in fact the set-wise stabilizer of $B\cap
\Gamma(C)$, and also of $B\cap \Gamma(D)$, in $X_B$. Then $K$ is
transitive on $B\cap \Gamma(C)$ and on $B\cap \Gamma(D)$. Let $H$
and $H_1$ be the kernels of $K$ acting on $B\cap \Gamma(C)$ and  on
$B\cap \Gamma(D)$, respectively. Then $K/H$ and $K/H_1$ are
permutation groups of degree $3$. Noting that $X_{(B)}\le H$ and
$X_{(B)}\le H_1$, it follows that $H/X_{(B)}$ and $H_1/X_{(B)}$ are
normal subgroups of $K/X_{(B)}$ with index $3$ in $K/X_{(B)}$. Hence
$H_1/X_{(B)}=H/X_{(B)}$ as $A_4$ has only one normal subgroup of
order $4$. Thus $H_1=H$ fixes $B$ point-wise, and so $H\le X_{(B)}$,
which contradicts $|H/X_{(B)}|=4$.

Suppose that $m^*(\Gamma,\mcB)\ne 1$. Recall that
$m^*(\Gamma,\mcB):=|B\cap(\cap_{C\in
\Gamma_\mathcal{B}(\mfv)}\Gamma(C))|$, the multiplicity of the dual
design $\mathcal{D}^*(B)$ of $\mathcal{D}(B)$, is independent of the
choices of $B$  and $\mfv\in B$. Assume that $\mathcal{D}^*(B)$ is a
$1$-$(v^*,b^*,r^*)$ design. Then $(v^*,b^*,r^*)=(b,v,k)$ is one of
$(4,4,3)$, $(4,6,3)$ and $(7,7,3)$. A similar argument as in the
above paragraph implies that $m^*(\Gamma,\mcB)$ is a divisor of $v$
and of $k$. Then $(b,v,k)=(4,6,3)$ and $m^*(\Gamma,\mcB)= 3=k$. It
follows that $m(\mcD(B))\ge |\Gamma_\mathcal{B}(\mfv)|=2$, again a
contradiction.

The above argument gives $m(\mcD(B))=1$ and $m^*(\Gamma,\mcB)=1$.
Then $X_B^{\Gamma_\mcB(B)}\cong X_B^B$  by
Theorem~\ref{thm_B_neighbour}. Thus $X_B^{\Gamma_\mcB(B)}$ is
$2$-transitive on $\Gamma_\mcB(B)$ if one of cases (a), (b) and (c)
occurs. Therefore, if one of Theorem \ref{thm_main}{(a)-(d)} occurs,
then $X_B$ acts $2$-transitively on the blocks of $\mathcal
{D}(B)=\Gamma_\mcB(B)$, and hence $\Gamma_{\mathcal {B}}$ is
$(X,2)$-arc-transitive.

\vskip 10pt

Now assume that $\Gamma_{\mathcal {B}}$ is $(X,2)$-arc-transitive.
Recall that $m(\mcD(B))$ is the multiplicity of $\mcD(B)$, the
number $C\in \Gamma_\mcB(B)$ with the same trace, which is
independent of the choice of $B$. Then $m(\mcD(B))=1$ by~\cite[Lemma
2.4]{FSGTTQ(II)}. Since $\Gamma_{\mathcal {B}}$ is
$(X,2)$-arc-transitive, $\lambda :=|\Gamma(C)\cap B\cap \Gamma(D)|$
is independent of the choice of $[C,B,D] \in Path_2(\Gamma_{\mathcal
{B}})$. By \cite[Corollary 3.3]{FSGTTQ(II)},   $vr = 3b$ and
$\lambda(b-1) = 3(r - 1)$, thus $(9 - \lambda v)r = 3(3-\lambda)$.
Since $\Gamma$ is not a multicover of $\Gamma_\mcB$, we have
$\lambda \leq k-1=2$ and $v>k$. If $\lambda = 0$, then $r = 1$ and
$v=3b$. Let $\lambda\ge 1$. Then, by \cite[Theorem 3.2]{FSGTTQ(II)},
the dual design $\mathcal{D}^*(B)$ of $\mathcal{D}(B)$ is a
$2$-$(b,r,\lambda)$ design with $v$ blocks. The well-known Fisher's
Inequality applied to $\mathcal{D}^*(B)$ gives $b\le v$, and so
$r\le k=3$. If $\lambda = 2$, then $\lambda(b-1) = 3(r - 1)$, $(9 -
2v)r = 3$ and $v>k$ imply   $(v,b, r) = (4,4, 3)$. If $\lambda = 1$,
then $ r \leq k $, $vr = 3b$ and $(9 - v)r = 6$ yield $(v,b,
r)=(6,4, 2)$ or $(7,7 ,3)$.

Note that $1\le m^*(\Gamma,\mcB)\le \lambda$ if $\lambda\ne
 0$. Suppose that $m^*(\Gamma,\mcB)\ne 1$ for some $\lambda\ne 0$.
 Then $\lambda=2=m^*(\Gamma,\mcB)$. It follows from $r=3=k$ that
 there are $C,D\in \Gamma_\mcB(\mfv)$ such that $C\ne D$ and $B\cap \Gamma(C)=B\cap
 \Gamma(D)$. Thus $C$ and $D$ has the same trace, and hence $m(\mcD(B))\ge 2$, a contradiction.
Therefore, if $\lambda\ne 0$ then $m^*(\Gamma,\mcB)=1$ and, by
Theorem~\ref{thm_quotient_star_design} and \ref{thm_B_neighbour},
$X_B^{\Gamma_\mcB(B)}\cong X_B^B$, and the induced action of $X$ on
$\mathcal {B}$ is faithful.

We treat four separate cases in the following.

{\bf Case 1}. $(v,b, r,\lambda) = (4,4, 3,2)$. Then
$val(\Gamma_{\mathcal {B}}) = 4$, and $X_B^B \cong A_4$ or $S_4$ as
$X_{B}$ acts $2$-transitively on $\Gamma_{\mathcal {B}}(B)$. Thus
(a) holds.

By \cite[Theorem 2]{ACOFSGWTTQ}, $\Gamma \cong \Xi(\Gamma_{\mathcal
{B}}, \Delta)$ for some self-paired $X$-orbit $\Delta$ on
$Arc_3(\Gamma_{\mathcal {B}})$. For any connected tetravalent
$(X,2)$-arc-transitive graph $\Sigma$, by Theorem
\ref{thm_evenvalency_selfpaired}, there exists some self-paired
$X$-orbit on $Arc_3(\Sigma)$, and by \cite[Theorem 10]{ACOFSGWTTQ},
the corresponding $3$-arc graph admits an $X$-invariant partition
with quotient graph isomorphic to $\Sigma$ and parameters
$(v,b,k,r)=(4,4,3,3)$. Thus, by Theorem \ref{thm_tetravalent}, one
of (a.1), (a.2) and (a.3) of Theorem \ref{thm_main}
 holds.

{\bf Case 2}. $(v,b, r,\lambda)=(6,4, 2,1)$. Then
$val(\Gamma_{\mathcal {B}}) = 4$, $X_B^B \cong A_4$ or $S_4$, and so
(b) occurs.

Since $(r, \lambda) = (2, 1)$, by Lemma \ref{lem_U-graph}, $\Gamma
\cong \gimel(\Gamma_{\mathcal {B}}, \Delta)$ for some self-paired
$X$-orbit $\Delta$ on $Arc_3(\Gamma_{\mathcal {B}})$. By Theorem
\ref{thm_tetravalent_self-paired}, $J(\Delta)$ is a self-paired
$X$-orbit on $J(\Gamma_{\mathcal {B}})$, and
$\gimel(\Gamma_{\mathcal {B}}, \Delta) \cong \Psi(\Gamma_{\mathcal
{B}}, J(\Delta))$. For any connected tetravalent
$(X,2)$-arc-transitive graph $\Sigma$,
 by Theorem \ref{thm_evenvalency_selfpaired}, there exists some self-paired
$X$-orbit on $Arc_3(\Sigma)$, and by Proposition~\ref{prop_U-graph}
and Theorem~\ref{thm_tetravalent}, the corresponding graph
constructed as in Proposition~\ref{prop_U-graph} admits an
$X$-invariant partition with quotient graph isomorphic to $\Sigma$
and parameters $(v,b, r,\lambda)=(6,4, 2,1)$. Then (b.1), (b.2) or
(b.3) follows from  Theorem \ref{thm_tetravalent},

{\bf Case 3}. $(v,b, r,\lambda)=(7,7 ,3,1)$. In this case,
$\mcD(B)\cong PG(2,2)$ is $X_B$-flag-transitive. Then
$X_B^{\Gamma_\mcB(B)}$ is isomorphic to a subgroup of $PSL(3,2)$,
the automorphism group of $PG(2,2)$. Since $\Gamma_{\mathcal {B}}$
is $(X,2)$-arc-transitive, $X_B^{\Gamma_\mcB(B)}$  is $2$-transitive
on $\Gamma_{\mathcal {B}}(B)$, and hence $|X_B^{\Gamma_\mcB(B)}|
\geq 42$. It follows that $X_B^{\Gamma_\mcB(B)}\cong PSL(3,2)$. Thus
$X_B^B \cong X_B^{\Gamma_{\mathcal {B}}(B)} \cong PSL(3,2)$  by
Theorem \ref{thm_B_neighbour}. Hence (c) holds. Since $ m^*
(\Gamma,\mcB)= 1$,  by Theorem \ref{thm_quotient_star_design},
$\Gamma \cong \Pi(\Gamma_{\mathcal {B}}, \Theta)$ for some
self-paired $X$-symmetric orbit $\Theta$ on $D\mcS
t^3(\Gamma_{\mathcal {B}})$. Further, by
Theorem~\ref{thm_heptavalent} and the above argument, one connected
heptavalent $(X,2)$-arc-transitive graph $\Sigma$ occurs as
$\Gamma_\mcB$ if and only if $X_\tau^{\Sigma(\tau)}\cong PSL(3,2)$.
Again by Theorem~\ref{thm_heptavalent}, one of (c.1), (c.2) and
(c.3) holds.

{\bf Case 4}. $\lambda=0$, $r=1$ and $v=3b$. Since $\Gamma_{\mathcal
{B}}$ is $(X,2)$-arc-transitive,  $X_B$ acts $2$-transitively on the
blocks of $\mathcal {D}(B)$. It follows from $r=1$ and $\lambda = 0$
that $\Gamma \cong e \Gamma[B, C]$ for $\{B, C\} \in
E(\Gamma_{\mathcal {B}})$. Since $k = 3$, we know $\Gamma[B,C]\cong
3K_2, K_{3,3}-3K_2$ or $K_{3,3}$, so one of   (d.1), (d.2) and (d.3)
occurs.


\vskip 10pt

\begin{thebibliography}{10}

\bibitem{AGT}
N.L. Biggs, Algebraic graph theory, Second Edition, Cambridge
Mathematical Library, Cambridge University Press, Cambridge, 1993.

\bibitem{DT} T. Beth, D. Jungnickel and H. Lenz, Design theory, Second Edition,
Cambridge University Press, 1999.

\bibitem{GTWA}
J.A. Bondy and U.S.R. Murty, Graph theory with applications, The
MaCmillan Press, London and Basingstoke, 1976.

\bibitem{PG}
J.D. Dixon and B. Mortimer, Permutation groups, Springer, New York,
1996.

\bibitem{FP1} X.G. Fang and C.E. Praeger, Finite two-arc transitive graphs
admitting a Suzuki simple group, {\em Comm. Algebra} {\bf 27}(1999),
3755-3769.

\bibitem{FP2}  X.G. Fang and C.E. Praeger, Finite two-arc transitive
graphs admitting a Ree simple group, {\em Comm. Algebra} {\bf
27}(1999), 3727-3754

\bibitem{Gardiner-Praeger95}
A.~Gardiner and C.E.~Praeger, A geometrical approach to imprimitive
graphs, {\em Proc. London Math. Soc.} (3) {\bf 71}(1995), 524-546.

\bibitem{Gardiner-Praeger98} A. Gardiner and C.E. Praeger, Topological covers of complete
graphs, {\em Math. Proc. Cambridge Philos. Soc.} {\bf 123}(1998),
549-559.

\bibitem{SGWCQ}
A.Gardiner and C.E.Praeger, Symmetric graphs with complete
quotients, preprint.


\bibitem{FSGWTTQ}
M.A. Iranmanesh, C.E. Praeger and S. Zhou, Finite symmetric graphs
with two-arc-transitive quotients, {\em J. Combin Theory B} {\bf
94}(2005), 79-99.


\bibitem{FSTGOOO}
C.H. Li, Finite s-arc-transitive graphs of odd order, {\em J.
Combin. Theory Ser. B} {\bf 81}(2001), 307-317.


\bibitem{ACOFSGWTTQ}
C.H. Li, C.E. Praeger and S. Zhou, A class of finite symmetric
graphs with $2$-arc-transitive quotients, {\em Math. Proc. Cambridge
Phil. Soc.} {\bf 129}(2000), 19-34.

\bibitem{OPPGWSSATOG}
C.H. Li, Z.P. Lu and D. Maru\v{s}i\v{c}, On primitive permutation
groups with small suborbits and their orbital graphs, {\em J.
Algebra} {\bf 279}(2004), 749-770.

\bibitem{FSGTTQ(II)}
Z.P. Lu and  S. Zhou, Finite Symmetric Graphs with Two-Arc
Transitive Quotients II, {\em J. of Graph Theory} {\bf 56}(2007),
167-193.

\bibitem{NG}
M. Perkel, Near-polygonal graphs, {\em Ars Combin.} {\bf
26}(A)(1988), 149-170.


\bibitem{Zhou2002a}
S. Zhou, Imprimitive symmetric graphs, 3-arc graphs and 1-designs,
{\em Discrete Math.} {\bf 244}(2002), 521-537.

\bibitem{Zhou2002b}
S. Zhou, Constructing a class of symmetric graphs, {\em Eur. J.
Combin.} {\bf 23} (2002), 741-760.

\bibitem{Zhou2004}
 S. Zhou, Almost covers of
2-arc-transitive graphs, {\em Combinatorica} {\bf 24}(2004),
731-745.

\bibitem{Zhou2005} S. Zhou, A local analysis of imprimitive
symmetric graphs, {\em J Algebraic Combin.} {\bf 22}(2005), 435-449.

\bibitem{OACOFSG}
S. Zhou, On a class of finite symmetric graphs, {\em Eur. J.
Combina.} {\bf 29}(2008), 630-640.

\bibitem{Zhoudm} S. Zhou, Classification of a family of symmetric graphs with
complete quotients, {\em Discrete Math.}, accepted.
\end{thebibliography}
\end{document}